\documentstyle[12pt,amsfonts]{article}
\newtheorem{theorem}{Theorem}[section]
\newtheorem{proposition}[theorem]{Proposition}

\newtheorem{lemma}[theorem]{Lemma}
\newtheorem{corollary}[theorem]{Corollary}

\newtheorem{definition}[theorem]{Definition}

\newcommand{\prf}{{\it Proof:} }
\newcommand{\rmk}{{\it Remark:} }

\newcommand{\Pic}{{\rm Pic}}

\renewcommand{\div}{{\rm div}}

\newcommand{\Spec}{{\rm Spec \,}}

\newcommand{\cC}{{\cal C}}
\newcommand{\cD}{{\cal D}}

\newcommand{\cF}{{\cal F}}
\newcommand{\cG}{{\cal G}}
\newcommand{\cH}{{\cal H}}

\newcommand{\cK}{{\cal K}}

\newcommand{\cO}{{\cal O}}

\newcommand{\cV}{{\cal V}}
\newcommand{\cW}{{\cal W}}
\newcommand{\cX}{{\cal X}}

\newcommand{\cZ}{{\cal Z}}
\newcommand{\A}{{\Bbb A}}

\newcommand{\C}{{\Bbb C}}

\newcommand{\E}{{\Bbb E}}

\renewcommand{\P}{{\Bbb P}\hspace{.05em}}
\newcommand{\Q}{{\Bbb Q}}
\newcommand{\R}{{\Bbb R}}

\newcommand{\Z}{{\Bbb Z}}

\newcommand{\qed}{\hspace*{\fill}\hbox{$\square$}}

\newcommand{\lra}{\longrightarrow}

\begin{document}
\title{\bf Constructing Indecomposable Motivic Cohomology Classes on
Algebraic Surfaces}
\author{Stefan J. M\"uller-Stach}
\date{version October 1995}
\maketitle

\begin{abstract}{ We describe a method to construct indecomposable classes in
Bloch's higher
Chow group $CH^2(X,1)$ on algebraic surfaces over $\C$ via transcendental
methods and apply it to obtain examples on K3 surfaces and some surfaces
of general type.  }
\end{abstract}
\ \\ \\
\section{Introduction}

Let $X$ be a smooth algebraic surface over $\C$. It is well known that
$K_0(X) \otimes \Q \cong \oplus CH^p(X) \otimes \Q$
as a consequence of Grothendieck's
Riemann-Roch theorem. S. Bloch has generalized this to higher algebraic
K-theory in \cite{Bl1}, see also \cite{Lev}. For example one obtains
$$K_1(X) \otimes \Q \cong \bigoplus_p CH^p(X,1)\otimes \Q $$
The groups $CH^p(X,n)$ are called {\it higher Chow groups}. \\
The purpose of this paper is to give an explicit way to construct
classes in $CH^2(X,1)$ that are non-trivial modulo the image of the natural
map
$$\gamma: {\rm Pic}(X) \otimes \C^* \to CH^2(X,1)$$
and modulo torsion. We call such classes {\it indecomposable}.
Cycles in $CH^2(X,1)$ can be constructed
by finding curves $Z_i $ in $X$ with rational functions $f_i$ on them that
satisfy $\sum {\rm div}(f_i)=0$ on $X$. By general conjectures
(see section 7), the cokernel of $\gamma$ is expected to be a countable
group on any smooth surface.\\
We will construct indecomposable elements in $CH^2(X,1)$
on general quartic K3-surfaces that contain a line
and on some special quintic surfaces of general type.
Other examples over the complex numbers have been constructed by M. Nori
(unpublished) on abelian surfaces, by A. Collino in \cite{Col} on Jacobian
varieties, and by C. Voisin and C. Oliva
on K3 surfaces in the unpublished work \cite{Voi2}.
The first examples ever given were over
number fields, by A. Beilinson in \cite{Be1}, other examples
on products of modular curves are contained in \cite{Fla}, \cite{Mil}
and probably at other places. \\
Our method consists of deforming the
complex structure of the pair $(X,Z)$ and studying the variation of mixed
Hodge structures associated to the open complements. The ideas of this
technique in the case of ordinary Chow groups are similar to those
developed in \cite{BMS} and in \cite{Voi3}, but eventually go back
to the fundamental idea of Griffiths to show the non-triviality of
cohomology classes on the general member of a family of varieties by
showing that their derivatives are non-zero.
The advantage of our method is that it is not restricted to surfaces with
trivial canonical bundle and that it is
very simple to apply in situations where some geometry is known, in
particular if an explicit family or a degeneration to a singular
configuration can be written down.\\
\ \\
Let us now describe the contents of this paper. In chapter 2 we first recall
the definition of Bloch's higher Chow groups
$CH^p(X,n)$ from \cite{Bl1} together with some of their basic properties.
Then we sketch the construction
of Deligne-Beilinson cohomology (\cite{Be1},\cite{EV})
and compare various definitions for the Chern class maps
$$c_{p,n}: CH^p(X,n) \to H^{2p-n}_\cD(X_{\rm an},\Z(p))  $$
($X_{\rm an}$ denotes the underlying analytic space, where we assume that $X$
is
defined over $\C$) due to Beilinson, Bloch, Gillet and others. In
particular we recall the explicit integral that computes $c_{2,1}$ and
study the relation between the extension of mixed Hodge structures given
by the Gysin sequence attached to the support of a given cycle
$Z \in CH^p(X,1)$ and $c_{p,1}(Z)$, using the description of Deligne
cohomology of smooth, projective varieties as ${\rm Ext^1}$ in the category
of mixed Hodge structures. \\
In chapter 3 we present the circle of ideas around the deformation
theory of Deligne-Beilinson cohomology classes and the rigidity of
Chern classes from higher Chow groups. This is essentially due to Bloch and
Beilinson, see \cite{Be1}. It implies in particular that
$$c_{p,n}: CH^p(X,n) \otimes \Q \to H^{2p-n}_\cD(X_{\rm an},\Q(p))  $$
has a countable image for $n \ge 2$ if $X$ is smooth and proper over $\C$.
In the case $n=1$ this is not true anymore, instead we can show the
following:
\begin{proposition} Let $X$ be a smooth and projective variety over $\C$.
Then for all $p$
the image of $c_{p,1}$ in $H^{2p-1}_\cD(X,\Q(p))/Hg^{p-1,p-1}(X) \otimes
\C/\Q(1)$ is countable.
\end{proposition}
Here $ Hg^{i,i}(X) \subset H^{2i}(X,\Z(i))$
denotes the set of Hodge classes.
This seems to be a well known fact, but I
could not  find a proof in the literature, so I decided to include one here
for the sake of completeness. A similar result holds in the case $n=0$,
stating that the Griffiths group has countable image in the intermediate
Jacobian $J^p(X)$ modulo the maximal abelian subvariety $J^p_a(X)$, see
\cite{Sh}. \\
\ \\
In chapter 4 these results are applied to the study of $CH^2(X,1)$ for an
algebraic surface $X$ over $\C$.
We prove the following result, which was communicated by H. Esnault:
\begin{proposition} Let $X$ be smooth, projective over $\C$. Then
$CH^2(X,1)$ decomposes iff ${\cal F}^2_\Z = {\cal F}^1_\Z \wedge
{\cal F}^1_\Z$ and $H^1(X,{\cal F}^2_\Z) \otimes \Q=0$.
\end{proposition}
The sheaves ${\cal F}^p_\Z$ have holomorphic p-forms with log-poles and
$\Z(p)-$periods as sections. A precise formulation can be found in the text.
The proof uses Gersten-Quillen resolutions. Then we discuss the relations
with Bloch's conjecture on the Chow groups of zero cycles on surfaces
mentioned already above and the result of Esnault-Levine (\cite{EL}) about
the relation between the decomposability of $CH^p(X,1)$ and the injectivity
of the cycle maps $c_{r,0}$ for $d-p+1 \le r \le d$, where $d={\rm dim}(X)$.
In the remaining part of the chapter we prepare the setup of variations
of mixed Hodge structures associated to a family of cycles $Z_t \in
CH^2(X_t,1)$.
Let us explain the necessary deformation theory. Assume we look at a
smooth, proper deformation $f:\cX \to S$ of $X$ with $S$ a smooth
and quasiprojective variety, a base point $0 \in S$ such that $f^{-1}(0)=X$
and a normal crossing divisor $\cZ$ in $\cX$  containing out of two smooth
components ${\cal Z}_1,{\cal Z}_2$ such that ${\cal Z}_1$ and
${\cal Z}_2$ (resp.
${\cal Z}_1 \cap {\cal Z}_2$) are smooth of relative dimension one (resp.
zero) over $S$ and restrict to $Z_1$ and $Z_2$ over the central fiber.
We get an exact diagram:
$$\matrix{ 0 & \to & T_X(log Z) & \to & T_{\cX}(log(\cZ))|_X
& \to & f^*T_{S,0} & \to & 0 \cr
& & \downarrow & & \downarrow  & & ||  & & \cr
0 & \to & T_X & \to & T_{ \cX}|X & \to & f^*T_{S,0} & \to & 0 }
$$
The {\it logarithmic
Kodaira-Spencer map}  is defined as the coboundary map
$$ T_{S,0} \longrightarrow H^1(X,T_X(log Z)) $$
Let us denote the image of $T_{S,0} \to H^1(X,T_X(log Z))$ by $W(log)$ and
the further image in $H^1(X,T_X)$ by $W$.\\
Our main result then is:

\begin{theorem} {\bf (CRITERION FOR INDECOMPOSABILITY)}: \\
Let $X$ be a smooth projective surface over $\C$.
Assume we are given two smooth and connected curves  $Z_1$ and $Z_2$
on $X$ intersecting
transversally and nontrivial rational functions $f_i$ on $Z_i$ ($i=1,2$),
such that ${\rm div}(f_1)+{\rm div}(f_2)=0$ as a zero-cycle on $X$. Denote by
$Z=Z_1 \otimes f_1 +Z_2 \otimes f_2$ the resulting cycle in
$CH^2(X,1)=H^1(X,\cK_2)$ and suppose the following
conditions hold:\\
(1) $Z$ also defines a cycle in Bloch's higher Chow group
$CH^1(|Z|,1)$ - again denoted by $Z$- and as such is not equivalent to
$Z_1 \otimes a_1 + Z_2 \otimes a_2$ with $a_1,a_2 \in \C^*$. \\
(2) There exist a smooth, proper deformation $f:\cX \to S$ with $S$ a smooth
and quasiprojective variety, a base point $0 \in S$ such that $f^{-1}(0)=X$
and the following properties hold:\\
(a) The situation in (1) deforms together with $X$: There exists a normal
crossing divisor $\cZ=\cZ_1+\cZ_2 \subset \cX$ with ${\cal Z}|_X=Z_1+Z_2$,
consisting out of two smooth components ${\cal Z}_1,{\cal Z}_2$
such that ${\cal Z}_1$ and ${\cal Z}_2$ (resp.
${\cal Z}_1 \cap {\cal Z}_2$) are smooth of relative dimension one (resp.
zero) over $S$. Furthermore there exist rational functions $F_i$ on $\cZ_i$
such that their restriction to each fiber $X_t:=f^{-1}(t)$ satisfy
${\rm div}(F_{1,t})+ {\rm div}(F_{2,t})=0$ as a zero-cycle in $X_t$ and
therefore define classes
$ Z_t=Z_{1,t} \otimes F_{1,t}+Z_{2,t} \otimes F_{2,t}$
in $CH^2(X_t,1)$ and in $CH^1(|Z_t|,1)$ for all $t \in S$.\\
(b)  The cup-product map
$$H^0(X,\Omega^2_X(logZ)) \otimes H^1(X,T_X(-Z)) \to H^1(X,\Omega^1_X)
/\oplus_i H^0(Z_i,\cO_{Z_i}) $$
has no left kernel.\\
(c) If $W(log) \subset H^1(X,T_X(logZ))$ denotes the image of
the logarithmic Kodaira-Spencer map in
$H^1(X,T_X(logZ))$,
then $W(log)$ contains the image of the natural map
$$H^1(X,T_X(-Z)) \to H^1(X,T_X(logZ))$$ \\
(d) For $t$ outside a countable number of proper analytic subsets of $S$,
$Z_{1,t}$ and $Z_{2,t}$ generate ${\rm NS}(X_t) \otimes \Q$.\\
{\bf Then:} $Z_t$ is non-torsion in $CH^2(X_t,1)/\Pic(X_t) \otimes \C^*$ for
$t$ outside a countable number of proper analytic subsets of $S$.
\end{theorem}
\ \\
Note that (2b) and (2d) together imply that the Picard number of $X_t$ is
not maximal for $t$ outside a countable number of proper analytic subsets of
$S$. Instead of the assumptions (2b) and (2c) we can also state a weaker
assumption (3) to get a sharper result:
\ \\
{\bf VARIANT:} \\
{\it Assume (1),(2a),(2d) of the above main theorem and additionally the
following instead of (2b),(2c):\\
(3)If $W(log) \subset H^1(X,T_X(logZ))$ denotes the image of
the logarithmic Kodaira-Spencer map in
$H^1(X,T_X(logZ))$, then the following map has no left kernel:
$$H^0(X,\Omega^2_X(logZ)) \otimes W(log) \to H^1(X,\Omega^1_X(logZ))
$$
{\bf Then:} $Z_t$ is non-torsion in $CH^2(X_t,1)/\Pic(X_t) \otimes \C^*$ for
$t$ outside a countable number of proper analytic subsets of $S$.}\\
\ \\
We prove both statements in chapter 5.
In chapter 6 we apply this result and Nori's connectivity theorem to study
some examples:\\
{\bf Example 1:}
\begin{theorem}
Let $X \subset \P^3$ be a general hypersurface of degree
$d \ge 5$. Then the Chern class map $CH^2(X,1) \otimes \Q \to
H^3_\cD(X,\Q(2)) $ has image isomorphic to $H^3_\cD(\P^3,\Q(2))
\cong \C/\Q(1)$.
\end{theorem}
This result should be seen as a generalization of the classical
Noether-Lefschetz theorem:
\begin{theorem}(Noether-Lefschetz) \\
Let $X \subset \P^3$ be a general hypersurface of degree
$d \ge 4$. Then $CH^1(X) \otimes \Q $ is isomorphic to
$H^2_\cD(\P^3,\Q(1)) \cong \Q $.
\end{theorem}
Both results can be shown using Nori's connectivity theorem \cite{Nori}
by the methods
of \cite{GM}. This was also observed by S.Bloch, M.Nori and C.Voisin
for example in \cite{Voi1}. We mention a slightly more general statement:
\begin{theorem}
Let $(Y,\cO(1))$ be a smooth and projective polarized variety
of dimension $n+h$, $X \subset Y$ a
general complete intersection of dimension $n$ and
multidegree $(d_1,..,d_h)$ with $min(d_i)$ sufficiently large. Furthermore
assume that $1 \le p \le n$. Then:
$${\rm Image}( CH^p(X,1) \otimes \Q \to H^{2p-1}_\cD(X,\Q(p))) $$
$$ \subset {\rm Image}(H^{2p-1}_\cD(Y,\Q(p))  \to H^{2p-1}_\cD(X,\Q(p))) $$
\end{theorem}

This result - already in the case of projective space - is somehow of
negative nature, because it destroys some obvious conjectures about the
image of higher Chern classes, see \cite{Voi1}. It leaves open the possibility
to construct examples on quartic K3-surfaces, also done in the paper
\cite{Voi2}. \\
{\bf Example 2:} We can show that the general quartic K3 surface $X$
containing a line has indecomposable $CH^2(X,1)$: Let $N$ be the
irreducible component
of smooth quartic surfaces containg a line such that some point $0 \in N$
corresponds to the Fermat quartic surface with the equation
$X_0=\{x \mid x_0^4+x_1^4-x_2^4-x_3^4=0 \} $.
Our cycles are obtained as follows:\\
A quartic surface $X$ that contains a line $G$ also contains an elliptic pencil
cut out by the residual elliptic curves of all hyperplane sections through
$G$. For a finite number of elements $E$ in this pencil, 2 of the 3
intersections points $P_1,P_2,P_3$ of $E$ and $G$ have the property that
$2(P_1-P_2)$ is rationally equivalent to zero. We show this in one example by
giving the explicit hyperelliptic map from $E$ to $\P^1$
ramified at those two points. Hence there is a
rational function $f_{1}$ on $E$ with zero divisor $2(P_1-P_2)$ and a rational
function $f_2$ on $G$ with zero divisor $2(P_2-P_1)$.
This construction can be extended to an irreducible
component $N$ of the Noether-Lefschetz locus of surfaces containing a
line. We obtain a cycle $E_t \otimes f_{1,t}+G_t \otimes f_{2,t} \in
CH^2(X_t,1)$ on a suitable covering $S$ of that component that respects the
choices of ordering of the chosen points.

\begin{theorem} Let $X_t$ be a general member of this family.
Then $CH^2(X_t,1)$ is not decomposable.
\end{theorem}
It is even true that the
general quartic K3-surface has indecomposable $CH^2(X,1)$, and this can be
proved by the same method as in the next example. The cycles used there were
used also by C. Oliva and C. Voisin in \cite{Voi2} on quartics. \\
In order to verify the assumptions of the main theorem we make use of
the Green- Gotzmann theorem and Griffiths' description of cohomology
groups of hypersurfaces via residues of differential forms, as
described in \cite{1594}. Assumption (2d) will hold for $t$ general on
such a component by Noether-Lefschetz theory.\\
With some more work, using a monodromy argument of H. Clemens which
was also used by A. Collino in \cite{Col}, one can probably prove infinite
generation for $CH^2(X_t,1)$ of a general quartic hypersurface
containing a line. Since this was also proved in \cite{Col} and \cite{Voi2}
we refrain from presenting it here. Obviously the idea would be to study the
monodromy around a countable set of loci on the parameter space
of the surfaces $X_t$.
\newpage
{\bf Example 3:} One can even obtain some examples of general
type: Look at the Shioda hypersurface of degree 5:
$$ X=\{x \in \P^3 \mid x_0x_1^4+x_1x_2^4+x_2x_0^4+x_3^5=0 \} $$
It has an automorphism $\sigma$ of order 65, given by
$$ \sigma: (x_0:x_1:x_2:x_3) \mapsto
(\zeta^{16}x_0:\zeta^{-4}x_1:\zeta x_2:x_3)$$
where $\zeta$ is a 65-th root of unity.
Shioda proves that the Picard group of $X$ is of rank one. Let $Z_1:=
X \cap H_3$ and $Z_2= X \cap H_0$,
where $H_i$ are the linear hyperplane sections $H_i=\{x_i=0\}$.
Then $Z_1 \cap Z_2$ intersect in two points (one
with multiplicity 4), called $P$ and $Q$ and one can show that $52P$ and
$52Q$ are rationally equivalent on both curves. $Z_2$ is not smooth, but
we construct a deformation $(X_t)_{t\in \C}$ and curves
$Z_{i,t}$ that are smooth for $t \neq 0$.
Taking a maximal irreducible
component $N$ of quintic surfaces such that it contains all $X_t$
and the cycle
$Z$ deforms along a suitable covering $S$ of $N$, the assumptions
of the theorem can be checked in the same way as in the previous example.
We obtain therefore:
\begin{theorem} A general member $X_t$ of the irreducible components of
quintics
that deform to the Shioda hypersurface $X$ and that preserve the given
cycle $Z \in CH^2(X,1)$ has indecomposable $CH^2(X_t,1)$.
\end{theorem}
In the remaining chapter 7 we sketch some ideas around these problems and
formulate some open problems. In particular we think that it would be
very convenient to have a good theory for singular surfaces in order to
get shorter proofs for indecomposability by degeneration methods.

\section{Higher Chow Groups and Chern Classes }

\subsection{Bloch's Higher Chow Groups }

Let $X$ be a quasiprojective variety over a field $k$. Define
$$\Delta^n := \Spec(k[T_0,...,T_n]/\sum T_i=1 ) $$
Then $\Delta^n \cong \A^n_k$ is affine n-space and by setting the coordinates
$t_i=1$ one obtains $(n+1)$ linear hypersurfaces in $\Delta^n$ called
codimension one faces. By iterating this one gets codimension $(n-m)$-faces
isomorphic to $\Delta^m$ for every $m<n$ inside of $\Delta^n$. These are
parametrized by strictly increasing maps $\rho:\{1,...,m\} \to \{1,...,n\}$.
Higher Chow groups are defined as the homology groups of a chain complex. Let
$Z^p(X,n) \subset Z^p(X \times \Delta^n)$ the subset of cycles
of codimension $p$ that meet all faces $X \times \Delta^m$ again in
codimension $p$ for $m<n$. Let $\partial_i: Z^p(X,n) \to Z^p(X,n-1)$ be the
restriction map to the $i-$th codimension one face for $i=0,...,n$ and let
$\partial=\sum (-1)^i \partial_i $. Then the homology of the complex
$$...\to Z^p(X,n+1) \to Z^p(X,n) \to Z^p(X,n-1) \to ...$$
at position $n$ is denoted by $CH^p(X,n)$ \cite{Bl1}. We will need the
following facts about higher Chow groups:  \\
(1) There is also a cubical version: Here let $\square^n:=(\P^1 \setminus
\{1\})^n $ with coordinates $t_i$ and codimension one faces
obtained by setting $t_i=0,\infty$. The rest of the definition is completely
analogous except that one has to divide out degenerate cycles
and it is known \cite{Bl1} that both complexes are quasiisomorphic.\\
(2) The groups $CH^*(X,*)$ are covariant for proper maps and contravariant
for flat maps.\\
(3) If $W \subset X$ is a codimension $r$ subvariety, then one has
localization
$$...\to CH^*(X,n) \to CH^*(X \setminus W,n) \to CH^{*-r}(W,n-r) \to
CH^*(X,n-1) \to ... $$    \\
(4) $CH^*(X,0)=CH^*(X)$ are the usual Chow groups.\\
(5) If $X$ is smooth, there is a product \cite{Bl1}
$$CH^p(X,q) \otimes CH^r(X,s) \to CH^{p+r}(X,q+s)$$
which can be easily defined using the cubical version. Thus it is
possible to define an action of correspondences on higher Chow groups.\\
(6) There exist cycle classes to Deligne-Beilinson cohomology \cite{Bl2}
in case $k$ is a field of characteristic zero: If we fix an embedding
$\sigma: k \hookrightarrow \C$ and denote by $X_{an}$
the associated complex analytic space then we have maps
$$c_{p,n}: CH^p(X,n) \to H^{2p-n}_\cD(X_{an},\Z(p))$$
They will be discussed in section (2.3).
\\
(7) There is a Riemann-Roch formula $K_n(X) \otimes \Q =\bigoplus_p CH^p(X,n)
\otimes \Q$, see \cite{Bl1} and \cite{Lev2}.\\
(8) Suslin's theorem \cite{Sus}: For $k$ itself: $CH^n(\Spec(k),n)=K_n^M(k)$
(Milnor K-theory).\\
(9) If $X$ is smooth and proper then $CH^1(X,1)=k^*$ and $CH^1(X,n)=0$
for $n \ge 2$. \\
(10) For $X$ smooth, we have $CH^p(X,1) =H^{p-1}(X,{\cK}_p)$ where
$\cK_p$ is Quillen's K-theory Zariski sheaf associated to the presheaf
$U \mapsto K_p(\cO(U)) $. Remember Bloch's formula $CH^p(X)=H^p(X,\cK_p)$.\\
\prf(for (10)) Consider the diagram
$$\matrix{ Z^p(X,2) & \to & Z^p(X,1) & \to & Z^p(X) \cr
           \downarrow N &   & \downarrow N & & || \cr
  {\oplus}_{x \in X^{(p-2)}} K_2(k(x)) & \to &
 \oplus_{x \in X^{(p-1)}} k(x)^* & \to & Z^p(X) }
$$

Here $N$ denotes the norm map. The inverse map is given by taking the graph
of rational functions on codimension $(p-1)-$subvarieties. To show that
the maps are mutually inverse to each other one needs only to show that
the graph of the norm of a cycle in $X \times \A^1$ is equivalent to
the cycle modulo $Z^2(X,2)$. This can be shown explicitely: First
reduce to the case where $X$ is a point and then use the explicit
formulas in \cite{Sus}.\qed

\rmk  By similar methods $CH^p(X,2) \to H^{p-2}(X,\cK_p)$ is surjective
and an isomorphism for $p=2$.

\subsection{Deligne-Beilinson Cohomology}
\label{Hg}
Let $X$ be any scheme of finite type over $\C$. Then there exists a twisted
duality theory in the sense of Bloch-Ogus consisting of Deligne homology
and cohomology groups (\cite{EV},\cite{Gil},\cite{Ja2})
$$H^i_\cD(X,\Z(j)),\quad H^\cD_i(X,\Z(j))  $$
satisfying the axioms in \cite{BO}. We just mention the duality isomorphism
$$ H^i_{\cD,Z}(X,\Z(j)) \cong H^\cD_{2d-i}(Z,\Z(d-j))$$
($d=\dim(X)$) for $X$ smooth and $Z \subset X$ a closed subvariety.
Very important for our purposes is also the weak purity statement under the
same assumptions:
The groups $H^i_{\cD,Z}(X,\Z(j))$ vanish for $i < 2r$ where $r$ denotes the
codimension of $Z$. Recall that for $X$ smooth one has an exact sequence
$$0 \to {H^{i-1}(X,\C) \over {F^j + H^{i-1}(X,\Z)}} \to
H^i_\cD(X,\Z(j)) \to F^j \cap H^i(X,\Z) \to 0 $$
where the notation $F^j \cap H^i(X,\Z)$ denotes the set of all classes
$\alpha \in H^i(X,\Z)$ such that $\alpha \otimes \C \in F^jH^i(X,\C)$. If
$X$ is smooth and proper we have additionally for $i=2j$:
$$0 \to J^j(X) \to H^{2j}_\cD(X,\Z(j)) \to F^j \cap H^{2j}(X,\Z) \to 0 $$
where $J^j(X)$ is the intermediate Jacobian. For $i < 2j$ and $X$ smooth
and proper, the group $F^j \cap H^i(X,\Z)$ is equal to the torsion
subgroup of $ H^i(X,\Z)$. By the isomorphism
$${\rm Ext}^1_{\rm MHS}(\Z(-j),H^{i-1}(X,\Z))=
{H^{i-1}(X,\C) \over {F^j + H^{i-1}(X,\Z)}}$$
both of the statements above can be subsumed into the exactness of the
sequence
$$0 \to {\rm Ext}^1_{\rm MHS}(\Z(-j),H^{i-1}(X,\Z)) \to H^i_\cD(X,\Z(j)) \to
{\rm Hom}_{\rm MHS}(\Z(-j),H^i(X,\Z)) \to 0 $$
for $X$ smooth and proper. See \cite{Be2}.

\subsection{Cycle Classes}

Let $k$ be a field of characteristic zero. Fix some embedding $\sigma: k
\hookrightarrow \C$ and let $X_{an}$ be the associated complex analytic
space. We give several definitions of cycle classes
$$c_{p,n}: CH^p(X,n) \to H^{2p-n}_\cD(X_{an},\Z(p))$$
which are in fact equivalent. There is Bloch's definition

\begin{definition} (\cite{Bl2})
\end{definition}
This is somehow the most general definition since it only uses some
functoriality and weak purity of Deligne cohomology
and can also be applied to get cycle classes
to \'etale cohomology. $X$ is just assumed to be quasiprojective over $k$.
A definition for all $c_{p,n}$ was given first by Beilinson
\begin{definition} (\cite{Be1},\cite{Gil},see \cite{Sch} for a survey)
\end{definition}
Again $X$ is quasiprojective over $k$. The definition uses simplicial schemes
and the axioms of Bloch and Ogus. The following definition is due to
Deninger and Scholl

\begin{definition} (\cite{DS},\cite{Ja1})
\end{definition}
Let $X$ be smooth and projective over $k$.
Here for $n \ge 1$ a cycle class
$$c_{p,n}: CH^p(X,n) \otimes \Q \to
{\rm Ext}^1_{\Q-{\rm MHS}}(\Q(-p),H^{2p-n-1}(X,\Q))$$
is defined by giving an explicit extension.
\ \\
These definition coincide by \cite{Bl2} for the first two definitions
and \cite{DS} for the first and third definition. If $n=1$ we can give more
detailed descriptions of the cycle class (\cite{Bl3},\cite{Lev}):\\
Assume that $X$ is smooth and projective over $\C$ and that
$H^{2p-1}(X,\Z)$ is torsionfree (otherwise work over $\Q$). Then we
have by property (10) of higher Chow groups that $CH^p(X,1)=H^{p-1}(X,\cK_p)$
and a class in this group is given by $Z=\sum Z_i \otimes f_i$ where the
$Z_i$ are integral subschemes of codimension $(p-1)$ and $f_i$ are
rational functions on each $Z_i$ with $\sum \div(f_i)=0$ as a cycle on $X$.
Then write $\div(f_i)=\partial \gamma_i$ for real analytic $(2d-2p+1)-$
dimensional chains on $Z_i$. In fact define $\gamma_i=f_i^{-1}(u_i)$ with
$u_i$ the standard path on the real axis from $0$ to $\infty$.
Let $\gamma=\cup \gamma_i$. We have $\partial
\gamma =0$ by $\sum \div(f_i)=0$ and even more $Z$ defines a class in
$F^pH^{2p-1}_{|Z|}(X,\Z)$ and therefore $\gamma$ has zero cohomology
class in  $H^{2p-1}(X,\Z)$ (by torsion-freeness). Thus
$\gamma =\partial \Gamma$ for some
chain $\Gamma \subset X$ of dimension $2d-2p+2$. The rational maps
$f_i:Z_i \to  \P^1$ are such that $\gamma_i $ is the preimage
of a path connecting $0$ and $\infty$ on $\P^1$ and hence allows to choose a
branch of logarithm on $Z_i \setminus \gamma_i$. The functional
$$ \alpha \mapsto \sum_i \int_{Z_i -\gamma_i} \log(f_i) \cdot \alpha +
(2\pi \sqrt{-1}) \int_\Gamma \alpha  $$
defines a functional on $(2p-2)$ forms on $X$ and therefore via Poincar\'e
duality a class in $H^{2p-2}(X,\C)$. The choices we made let it become well
defined only in the quotient group $H^{2p-1}_\cD(X,\Z(p))$.
\ \\
The definition via extensions can also be described partially by the
following construction:\\
If we set $U=X \setminus \cup Z_i$, then one has a long exact sequence
$$  H^{2p-2}_{|Z|}(X,\Z) \to H^{2p-2}(X,\Z) \to H^{2p-2}(U,\Z)
\to H^{2p-1}_{|Z|}(X,\Z) \to H^{2p-1}(X,\Z) $$
A cycle $Z \in CH^p(X,1)$ gives rise to a class in $F^p \cap
H^{2p-1}_{|Z|}(X,\Z)$ and hence gives an extension
$$ 0 \to H^{2p-2}(X,\Z)/Hg^{p-1,p-1}(X) \to \E \to \Z(-p) \to 0 $$
Here $Hg^{i,i}=F^i \cap H^{2i}(X,\Z)$ with the meaning defined in \ref{Hg}, in
particular $Hg^{1,1}={\rm NS}(X)$.
This extension class is the image of $c_{p,1}(Z)$ under the map
$${\rm Ext}^1(\Z(-p),H^{2p-2}(X,\Z)) \to
{\rm Ext}^1(\Z(-p),H^{2p-2}(X,\Z)/Hg^{p-1,p-1}(X))$$

\rmk Once we believe in the existence of regulator maps also for the local
situation, there is another way to describe the cycle classes for $n=1$
(\cite{E2}):\\
Let $\cK_p$ the Quillen K-theory sheaf (see above) and $\cH_\cD^p(p)$ be the
sheafified Deligne-Beilinson cohomology
(with presheaf $U \mapsto H^p_\cD(U,\Z(p))$) both
viewed as sheaves in Zariski topology.
Then there is a Leray spectral sequence (as in \cite{BO}) arising
from changing from Zariski to the analytic site
$$E_2^{p,q}(r)=H^p_{\rm Zar}(X,\cH_\cD^q(r)) \Longrightarrow
H^{p+q}_{\cD,\rm an}(X,\Z(r))  $$
and we get an edge morphism
$H^{p-1}(X,\cH_\cD^p(p)) \to H^{2p-1}_\cD(X,\Z(p))$
as a byproduct. The existence of regulator maps on affine schemes implies
a regulator map of sheaves
$${\rm reg}: \cK_p \to \cH_\cD^p(p)$$
The composition
$CH^p(X,1) \to H^{p-1}(X,\cK_p) \to H^{p-1}(X,\cH_\cD^p(p))
\to H^{2p-1}_\cD(X,\Z(p))$
is also equivalent to the cycle classes defined above.\\

\section{Deformations and Rigidity of Cycle Classes}

Let $Y$ be a reduced quasiprojective scheme over $\C$ and $A$ a ring with
$\Q \subset A \subset \R$. Let $\epsilon_A$ be the map
$\epsilon_A: H^k_\cD(Y,A(p)) \to F^p \cap H^k(Y,A(p))$.

\begin{lemma} (\cite{Be1}, 1.6.6.1.)\\
(a) If $p > min(k,\dim(Y))$, then $\epsilon_A \equiv 0$.\\
(b) If $Y=X \times S$ with $X$ smooth, projective and $S$ smooth, affine and
$k < 2p-\dim(S)$ then also $\epsilon_A \equiv 0$.
\end{lemma}

\prf In both cases $F^p \cap H^k(Y,A(p))=0$ by type reasoning.\qed

\begin{lemma} (\cite{Be1}, 1.6.6.2.)\\
Let $X$ be projective, $S$ smooth affine, $Y=X \times S$ and given
$s_1,s_2 \in S$. Assume
$k \le 2p-2$ and a class $\alpha \in H^k_\cD(Y,A(p))$ is given. Then:
$$\alpha|_{X \times \{s_1\}} = \alpha|_{X \times \{s_2\}} \in
H^k_\cD(X,A(p))$$
\end{lemma}

\prf $S$ is an affine curve wlog. By the lemma above
$\epsilon_A(\alpha)=0$, d.h. $\alpha \in {H^{k-1}(Y,\C) \over {F^p
\oplus H^{k-1}(Y,A(p))}}$. But Betti classes are rigid. \qed

\begin{corollary} (\cite{Be1}, 2.3.4.)\\
Let $X$ be smooth, projective over $\C$ and  $n \ge 2$.
Then the image of
$c_{p,n}:CH^p(X,n) \otimes \Q \to H^{2p-n}_\cD(X,\Q(p))$ is countable.
\end{corollary}
\prf There exists an algebraically closed, countable field $L \subset \C$
such that $X = X_0 \otimes_L \C$ for some $L-$variety $X_0$. Hence
$CH^p(X_0,n)$ is countable and therefore the image of
$c_{p,n}: CH^p(X_0,n)\otimes \Q \to H^{2p-n}_\cD(X_{an},\Q(p)) $.
It remains to show that the cycle classes from $X$ and $X_0$
have the same image in $H^{2p-n}_\cD(X_{an},\Q(p)) $.
Choose a smooth affine scheme $S=Spec(R)$ for some $L-$algebra $R$ and a
L-rational point $0\in S$, such that  the geometric general fiber in
$X_0 \times S$
is given by $X_\C$ . Given a cycle $Z \in CH^p(X,n) \otimes \Q$,
there is a spreading
${\cal Z} \in CH^p(X_0 \times S,n) \otimes \Q$
with Deligne class
$\alpha \in H^{2p-n}_\cD(X_{an} \times S_{an},\Q(p))$. In particular this
means that $Z={\cal Z}|_X$, the restriction to the geometric general fiber.\\
By the lemma above $\alpha|_X=\alpha|_{X_0}$ and hence $c_{p,n}(Z)$ is also
contained in the countable image of
$c_{p,n}(CH^p(X_0,n) \otimes \Q) \subset H^{2p-n}_\cD(X_{an},\Q(p))$. \qed

Now consider the cases $k=2p-1,2p$. Let $X$ be again smooth and projective
over $\C$ and $S$ an affine, smooth complex curve with good
compactification ${\overline S}=S \cup \Sigma $. Then:\\

\begin{lemma}

(a) If $n=1$, then $$F^p \cap H^{2p-1}(X \times S,A(p)) =
    H^{p-1,p-1}(X,A(p-1)) \oplus \bigoplus_\Sigma A(-1)
    $$ where the summation ranges over all divisors at infinity.\\
(b) (see \cite{ES}.) If $n=0$, then  $$F^p \cap H^{2p}(X \times S,A(p)) =
    F^p \cap H^{2p}(X \times {\overline S},A(p))
    \subset  $$
$$ \subset H^{p,p}(X,A(p)) \oplus
    H^{p-1,p}(X,\C)\otimes H^0({\overline S},
    \Omega^1_{\overline S}(log \Sigma)) \oplus
    H^{p,p-1}(X,\C)\otimes H^1({\overline S},{\cal O}_{\overline S}) $$
\end{lemma}
\prf In both cases
$F^{p+1} \cap H^{k}(X \times S,A(p)) = 0$ and the weight $k$ piece arises
exactly from the compactification of $X \times S$. This is the only weight
in (b). The claim follows therefore from K\"unneth decomposition. (a) follows
in the same way but there is only the weight $k+1=2p$.  \qed

\rmk For
$k=2p-1$, $\epsilon(\alpha)$ is only determined by the residues around
boundary divisors.
\\
If $k=2p$, then $\epsilon(\alpha)=\alpha_1+\alpha_2+\alpha_3$,
with $\alpha_1 \in H^{p,p}(X,A(p))$ the cohomology class of a restriction
to the general fiber.
$\alpha_2$ and $\alpha_3$ are complex conjugate.\\

\begin{corollary}
Let $X$ be smooth, projective over $\C$ and  $n=1$.
Then the image of the truncated cycle class
$$c_{p,1}^{\rm tr}:CH^p(X,1) \otimes \Q \to
H^{2p-1}_\cD(X_{an},\Q(p))/{\rm Hg}^{p-1,p-1}(X) \otimes \C/\Q(1)$$
is countable.
\end{corollary}
\prf Let $S$ be a smooth, affine and connected curve.
Consider the product map
$$ CH^{p-1}(X) \otimes CH^1(S,1) \to CH^p(X \times S,1) $$
respectively its Deligne cohomology version
$$ \gamma: H^{2p-2}_\cD(X_{\rm an},\Z(p-1)) \otimes H^1_\cD(S_{\rm an},\Z(1))
\to H_\cD^{2p-1}(X_{\rm an} \times S_{\rm an},\Z(p))  $$
Therefore
$H_\cD^{2p-1}(X_{\rm an} \times S_{\rm an},\Z(p))/ {\rm Image}(\gamma)$
becomes equal to
$${ H^{2p-2}(X \times S,\C^*) \over
{F^p + {\rm Hg}^{p-1,p-1}(X) \otimes \C^*
+{\rm other \quad terms}}}$$

If we mod out by the image of $\gamma$, we get a restriction map (the other
terms restrict to zero) for every $t \in S$:
$$ r_t: {H_\cD^{2p-1}(X_{\rm an} \times S_{\rm an},\Z(p)) \over
{\rm Image}(\gamma) }
\to H^{2p-1}_\cD(X_{an},\Z(p))/{\rm Hg}^{p-1,p-1}(X) \otimes \C^*$$
with the important property that the truncated cycle class $c_{p,1}^{\rm tr}$
factors through it. Therefore
we again have rigidity: Given $\alpha \in
H_\cD^{2p-1}(X_{\rm an} \times S_{\rm an},\Z(p))/ {\rm Image}(\gamma)$
one has after tensoring with $\Q$:
$\alpha|_{X \times \{s_1\}} = \alpha|_{X \times \{s_2\}}$ in
$H^{2p-1}_\cD(X_{an},\Q(p))/{\rm Hg}^{p-1,p-1}(X) \otimes \C/\Q(1)$, since
$H_\cD^{2p-1}(X_{\rm an} \times S_{\rm an},\Z(p))/ {\rm Image}(\gamma)$
can be represented by Betti classes.
Now the same argument as in the proof of Cor. 3.3. can be applied to
deduce the countability.
\qed
\section{Theory of $CH^2(X,1)$}

\subsection{Decomposability}

Let $X$ be a smooth and projective variety over an algebraically closed field
of characteristic zero. There is a natural map
$$ \gamma: \Pic(X) \otimes k^* \lra CH^2(X,1)  $$
which in the cubical version of higher Chow groups can be described by
sending $D \otimes a $ to $D \times \{a\} \subset X \times \square^1 $ for
$D$ an integral subscheme and $a \in k$ and extending linearly. \\
Alternatively the cup-product map
$\cO_X^* \otimes_\Z \cO_X^* \to \cK_2 $ gives rise to the map
$$H^1(X,\cO_X^*) \otimes H^0(X,\cO_X^*) \to H^1(X,\cK_2) $$
Equally consider the product map
$$CH^1(X,0) \otimes CH^1(X,1) \to CH^2(X,1)  $$
and it is an easy exercise to show that all three definitions coincide.

\begin{definition} We say that $CH^2(X,1)$ is decomposable, if the map
$ \gamma $ has torsion cokernel. More generally
we will say that $CH^p(X,1) $ is decomposable if the map
$CH^{p-1}(X) \otimes CH^1(X,1) \to CH^p(X,1) $ has torsion cokernel.
\end{definition}
The cokernel is a birational invariant in that case.
About the kernel we note that the composed map
${\rm Pic}^0(X) \otimes k^* \to CH^2(X,1) \to H^3_\cD(X_{an},\Z(2)) $
is zero, however we do not know in which way the map
${\rm Pic}^0(X) \otimes k^* \to CH^2(X,1)$ itself behaves.

\subsection{Criteria for Decomposability}

Let $X$ be as in the previous section and $k=\C$.
Let us fix some notations:\\
Denote by $\cF^p_\Z$ the Zariski sheaf associated to the presheaf
that associates to each open set $U$ the vector space of holomorphic
$p-$forms with $\Z(p)-$ periods and logarithmic poles along a
desingularization of $X \setminus U$. One has an exact sequence
$$0 \to \cH^{p-1}_{\rm DR}(\C/\Z(p)) \to \cH^p_\cD(p) \to \cF^p_\Z \to 0 $$
Here $\cH^{p-1}_{\rm DR}(\C/\Z(p))$ is the sheaf associated to the
presheaf
$U \mapsto H^{p-1}(U,\C/\Z(p))$ and
$\cH^p_\cD(p)$ is the sheafified Deligne-Beilinson cohomology as explained
before. We then have $\cF^1_\Z={\rm dlog}\cO^*_X$ and define
$\cC:= \cF^2_\Z /\cF^1_\Z \wedge \cF^1_\Z $. One can show that the sheaves
$\cH^{p-1}_{\rm DR}(\C/\Z(p)), \cH^p_\cD(p)$ and $\cF^p_\Z$ admit
Gersten-Quillen type resolutions which we will use below.
The following is
private communication by H. Esnault and partially explained in \cite{E1}.

\begin{theorem}  $CH^2(X,1)$ decomposes if
$\cC=0$ and $H^1(X,\cF^2_\Z)\otimes \Q=0$.
\end{theorem}

\prf
Consider the surjective map
$\alpha: \cK_2 \to {\cal F}^1_\Z \wedge {\cal F}^1_\Z$ induced by the
dlog map. By assumption $H^1(X,{\cal F}^1_\Z \wedge {\cal F}^1_\Z) \otimes
\Q=0$ and therefore it is sufficient to show that $H^1(X,\cK_2^0)$
decomposes where $\cK_2^0:={\rm Ker}(\alpha)$.
Gersten-Quillen resolutions for all three sheaves form a
commutative diagram
$$\matrix{ \cK^2_0 & \to  & K_2^0(\C(X)) & \to &
\oplus_{X^{(1)}} \C/\Z(1) & & \cr
\downarrow && \downarrow && \downarrow &&  \cr
\cK_2 & \to & K_2(\C(X)) & \to &
\oplus_{X^{(1)}} K_1(\C(D)) & \to & \oplus_{X^{(2)}} \Z \cr
\downarrow && \downarrow && \downarrow && || \cr
\cF^2_\Z & \to & F^2_\Z(\C(X)) & \to &
\oplus_{X^{(1)}} F^1_\Z(\C(D)) & \to & \oplus_{X^{(2)}} \Z     }
$$
where the third column is the direct sum of exact sequences
$$0 \to \C/\Z(1) \to \C(D)^*  {\buildrel {\rm dlog} \over \longrightarrow}
F^1_\Z(\C(D)) \to 0 $$
Thus we have a surjective map ${\rm Pic(X)} \otimes \C^* \to
H^1(X,\cK_2^0)$
and the assertion follows.    \qed

\rmk The statement in the theorem holds also in the reverse direction.
Bloch conjectures that for a complex algebraic surface
with $p_g=0$ the kernel of the Albanese map $CH^2_0(X) \to {\rm Alb}(X)$
is zero. The conjecture holds if $X$ is not of general type by \cite{BKL}.
If $X$ is of general type $p_g=0$ implies also $q(X)=0$ and the conjecture
has been verified for Godeaux- and Barlow surfaces for example.
Note that if $p_g >0$, then by \cite{Mum} quite
the contrary happens. If
Bloch's conjecture holds, it implies that $CH^2(X,1)$ decomposes
if $p_g(X)=0$ by \cite{BS}.
There is a more general statement in \cite{EL}:
\begin{theorem} (\cite{EL}) \\
Let $X$ be a smooth algebraic variety over $\C$.
Suppose that the cycle maps
$$c_{p,0}: CH^p(X) \to H^{2p}_\cD(X,\Z(p))$$
are injective for $d-s \le p \le d=\dim(X)$ for some $s \ge 0$.
Then $CH^p(X,1)$ is decomposable for $0 \le p \le s+1$.
\end{theorem}

\subsection{Criteria for Non-Decomposability and Main Theorem }

In this section assume for simplicity that $X$ is an algebraic surface over
$\C$. To show that $CH^2(X,1)$ is not decomposable we use cycle class maps to
Deligne cohomology. Note that the image of
$$c_{1,2}:CH^2(X,1) \to H^3_\cD(X,\Z(2))  $$
when restricted to $\Pic(X) \otimes \C^*$ is contained in
${\rm NS}(X) \otimes \C^*$ considered as a subgroup of $H^3_\cD(X,\Z(2))$.
Therefore it will be enough to show that there are classes not contained in
that subgroup modulo torsion. In fact the image of
$CH^2(X,1)/\Pic(X) \otimes \C^*$ in
$H^3_\cD(X,\Z(2))/{\rm NS}(X) \otimes \C^*$
is at most countable by Corollary (3.5).\\
In the last section we will construct examples and the
following theorem provides the necessary technical tool. We use the following
notation: If $Z$ is a divisor with normal crossings on $X$ and
smooth components $Z_i$, then let
$\Omega^p_X(logZ)$ be the sheaf of holomorphic p-forms with log-poles along
$Z$ and $T_X(logZ)$ be the dual of $\Omega^1_X(logZ)$.
There is an exact sequence
$$0 \to \Omega^1_X \to \Omega^1_X(logZ) \to \oplus \cO_{Z_i} \to 0 $$
and a commutative diagram where we define
$\cG:=\Omega^2_X(logZ)/\Omega^2_X$:
$$\matrix{\Omega^2_X &\to &W_1 \Omega^2_X(logZ) &\to &\oplus \Omega^1_{Z_i}
\cr
|| & & \downarrow & & \downarrow \cr
\Omega^2_X & \to & \Omega^2_X(logZ) & \to & \cG \cr
& & \downarrow & & \downarrow \cr
& & \C^k & = & \C^k }
$$
where $\C^k$ is supported on the $k$ points in the singular locus of $Z$
consisting of the intersection points in $Z_i \cap Z_j$.
The cup product with the Kodaira-Spencer class
$H^0(X,\Omega^2_X) \otimes H^1(X,T_X) \to H^1(X,\Omega^1_X) $ and
the natural map
$H^1(X,T_X(logZ)) \to H^1(X,T_X)$ induce a commutative diagram
$$\matrix{
H^0(X,\Omega^2_X) \otimes H^1(X,T_X(logZ)) & \to & H^1(X,\Omega^1_X) \cr
\downarrow & & \downarrow \cr
H^0(X,\Omega^2_X(logZ)) \otimes H^1(X,T_X(logZ)) & \to &
H^1(X,\Omega^1_X(logZ)) }
$$
Let us explain the necessary deformation theory. Assume we look at a
smooth, proper deformation $f:\cX \to S$ of $X$ with $S$ a smooth
and quasiprojective variety, a base point $0 \in S$ such that $f^{-1}(0)=X$
and a normal crossing divisor $\cZ$ in $\cX$  containing out of two smooth
components ${\cal Z}_1,{\cal Z}_2$ such that ${\cal Z}_1$ and
${\cal Z}_2$ (resp. ${\cal Z}_1 \cap {\cal Z}_2$)
are smooth of relative dimension one (resp.
zero) over $S$ and restrict to $Z_1$ and $Z_2$ over the central fiber.
We get an exact diagram
$$\matrix{ 0 & \to & T_X(log Z) & \to & T_{\cX}(log(\cZ))|_X
& \to & f^*T_{S,0} & \to & 0 \cr
& & \downarrow & & \downarrow  & & || & & \cr
0 & \to & T_X & \to & T_{\cX}|X & \to & f^*T_{S,0} & \to & 0 }
$$
The {\it logarithmic
Kodaira-Spencer map}  is defined as the coboundary map
$$ T_{S,0} \longrightarrow H^1(X,T_X(log Z)) $$
Let us denote the image of $T_{S,0} \to H^1(X,T_X(log Z))$ by $W(log)$ and
the further image in $H^1(X,T_X)$ by $W$. \\
Before stating the theorem let  us desribe one of the higher Chow groups of
the projective (however reducible) variety $|Z|$, the
support of $Z$. For simplicity we will
assume that $|Z|$ consists of two smooth components $|Z|=Z_1+Z_2$ with $k$
intersection points in $Z_1 \cap Z_2$.

\begin{proposition}. $CH^1(|Z|,1) \cong H^3_{\cD,|Z|}(X,\Z(2))$ and there
is an exact sequence
$$0 \to (\C^*)^{\oplus 2} \to CH^1(|Z|,1)
{\buildrel \tau \over \to}
{\rm Ker}(\Z^k \to \oplus {\rm Pic}(Z_i))  \to 0 $$
where $\Z^k$ is supported on $Z_1 \cap Z_2$.
\end{proposition}
\prf next section. \qed

Let us denote the map $CH^1(|Z|,1) \to {\rm Ker}(\Z^{k} \to
\oplus {\rm Pic}(Z_i))$ given in the proposition by $\tau$. Now we are
ready to state the main result and a stronger variant of it:

{\bf MAIN THEOREM (CRITERION FOR INDECOMPOSABILITY):} \\
{\it Let $X$ be a smooth projective surface over $\C$.
Assume we are given two smooth and connected curves  $Z_1$ and $Z_2$
on $X$ intersecting
transversally and nontrivial rational functions $f_i$ on $Z_i$ ($i=1,2$),
such that ${\rm div}(f_1)+{\rm div}(f_2)=0$ as a zero-cycle on $X$. Denote by
$Z=Z_1 \otimes f_1 +Z_2 \otimes f_2$ the resulting cycle in
$CH^2(X,1)=H^1(X,\cK_2)$ and suppose the following
conditions hold:\\
(1) $Z$ also defines a cycle in Bloch's higher Chow group
$CH^1(|Z|,1)$ - again denoted by $Z$- and as such is not equivalent to
$Z_1 \otimes a_1 + Z_2 \otimes a_2$ with $a_1,a_2 \in \C^*$. \\
(2) There exist a smooth, proper deformation $f:\cX \to S$ with $S$ a smooth
and quasiprojective variety, a base point $0 \in S$ such that $f^{-1}(0)=X$
and the following properties hold:\\
(a) The situation in (1) deforms together with $X$: There exists a normal
crossing divisor $\cZ=\cZ_1+\cZ_2 \subset \cX$ with ${\cal Z}|_X=Z_1+Z_2$,
consisting out of two smooth components ${\cal Z}_1,{\cal Z}_2$
such that ${\cal Z}_1$ and ${\cal Z}_2$ (resp.
${\cal Z}_1 \cap {\cal Z}_2$) are smooth of relative dimension one (resp.
zero) over $S$. Furthermore there exist rational functions $F_i$ on $\cZ_i$
such that their restriction to each fiber $X_t:=f^{-1}(t)$ satisfy
${\rm div}(F_{1,t})+ {\rm div}(F_{2,t})=0$ as a zero-cycle in $X_t$ and
therefore define classes
$ Z_t=Z_{1,t} \otimes F_{1,t}+Z_{2,t} \otimes F_{2,t}$
in $CH^2(X_t,1)$ and in $CH^1(|Z_t|,1)$ for all $t \in S$.\\
(b)  The cup-product map
$$H^0(X,\Omega^2_X(logZ)) \otimes H^1(X,T_X(-Z)) \to H^1(X,\Omega^1_X)
/\oplus_i H^0(Z_i,\cO_{Z_i}) $$
has no left kernel.\\
(c) If $W(log) \subset H^1(X,T_X(logZ))$ denotes the image of
the logarithmic Kodaira-Spencer map in
$H^1(X,T_X(logZ))$,
then $W(log)$ contains the image of the natural map
$$H^1(X,T_X(-Z)) \to H^1(X,T_X(logZ))$$\\
(d) For $t$ outside a countable number of proper analytic subsets of $S$,
$Z_{1,t}$ and $Z_{2,t}$ generate ${\rm NS}(X_t) \otimes \Q$.\\
{\bf Then:} $Z_t$ is non-torsion in $CH^2(X_t,1)/\Pic(X_t) \otimes \C^*$ for
$t$ outside a countable number of proper analytic subsets of $S$.}
\ \\ \\
{\bf VARIANT:} \\
{\it Assume (1),(2a),(2d) of the above main theorem and additionally the
following instead of (2b),(2c):\\
(3)If $W(log) \subset H^1(X,T_X(logZ))$ denotes the image of
the logarithmic Kodaira-Spencer map in
$H^1(X,T_X(logZ))$, then the following map has no left kernel:
$$H^0(X,\Omega^2_X(logZ)) \otimes W(log) \to H^1(X,\Omega^1_X(logZ))
$$
{\bf Then:} $Z_t$ is non-torsion in $CH^2(X_t,1)/\Pic(X_t) \otimes \C^*$ for
$t$ outside a countable number of proper analytic subsets of $S$.}\\
\ \\ \\
\rmk The assumption that we have only two components is not necessary but
simplifies the proof and suffices for the  applications. (2b) and (2d)
imply necessarily that $p_g(X) \ge 1$, since it follows from
(2b) and (2d) that the Picard number of $X_t$ is not maximal for $t$ outside
a countable number of analytic subsets of $S$.
(2b) and (2c) imply (3) (see proof of the main theorem) and therefore the
variant is the more general formulation.

\section{Proof of the Main Theorem }

\subsection{Auxiliary Results in Hodge Theory}

Let $S$ be a smooth complex variety.

\begin{definition} (\cite{SZ})\\
A graded polarized $\Z-$variation of mixed Hodge structures ($\Z$-VMHS)
on $S$ is a local system $\cV$ on
$S$ of $\Z-$modules of finite rank with the following data:\\
(a) An increasing filtration ...$\subset \cW_k \subset \cW_{k+1} \subset
\cV \otimes_\Z \Q $ by local systems over $\Q$.  \\
(b) A decreasing filtration $...\subset \cF^{p+1} \subset \cF^p \subset
... \subset \cF^0=\cV \otimes_\Z \cO_S$ by holomorphic vector bundles.\\
(c) The flat connection $\nabla$ on $\cV$ satisfies $\nabla \cF^p \subset
\Omega^1_S \otimes \cF^{p-1}$.\\
(d) The $\cF-$filtration induces on the local systems ${\rm Gr}^W_k \cV_\Q:=
\cW_k/\cW_{k-1}$ a pure variation of polarized Hodge structures on $S$.
\end{definition}

The graded pieces of the connection we denote by $\nabla^p$:
$$\nabla^p: \cF^p/\cF^{p+1} \to \Omega^1_S \otimes \cF^{p-1}/\cF^p  $$
In particular there are examples where such VMHS come from a geometric
situation, for example in the situation of the theorem:
Assume we are given an explicit
smooth, proper deformation $f:\cX \to S$ with $S$ a smooth
complex variety and there exist a cycle
$\cZ \in CH^2(\cX,1)$ with ${\cal Z}|_X=Z$ consisting out of two
smooth components ${\cal Z}_1,{\cal Z}_2$ intersecting transversally,
such that ${\cal Z}_1$ and ${\cal Z}_2$ (resp.
${\cal Z}_1 \cap {\cal Z}_2$) are smooth of relative dimension one (resp.
zero) over $S$ and such that its restriction to each fiber
$X_t:=f^{-1}(t)$ defines classes $ Z_t \in H^1(X_t,\cK_2)$.
Then it is a result in \cite{SZ} that
the cohomology groups $H^2(X_t \setminus Z_t,\Z)$ form an
admissible $\Z$-VMHS over $S$. \\
If all fibers are algebraic surfaces, there are no holomorphic 3-forms, hence
$\cF^3=0$ and $\cF^2$
is the holomorphic subbundle with fibers $H^0(X_t,\Omega^2_{X_t}(logZ_t))$.
The graded piece $\nabla^2: \cF^2 \to \Omega^1_S
\otimes \cF^1/\cF^2$ can be described in the stalk at $0 \in S$ by the maps
$$
\nabla^2: H^0(X,\Omega^2_X(logZ)) \to  H^1(X,\Omega^1_X(logZ)) \otimes
W(log)^*
$$
where $W(log) $ is the logarithmic Kodaira-Spencer image of $T_{S,0}$.
Bringing $W(log)$ to the other side we get the cup-product maps
$$
H^0(X,\Omega^2_X(logZ)) \otimes W(log) \to  H^1(X,\Omega^1_X(logZ))
$$
\begin{lemma} In this situation, assume that \\
$H^0(X,\Omega^2_X(logZ)) \otimes W(log) \to  H^1(X,\Omega^1_X(logZ)) $ has
no left kernel.
\\
{\bf Then:} $F^2 \cap H^2(U_t,\Q)=0$ for $t$ outside a countable number
of analytic subsets of $S$.
\end{lemma}
\prf  Assume there exists a nonzero
class $\lambda_0 \in F^2 \cap H^2(U_0,\Q)$ that
is supported over some germ of an analytic subvariety $S(\lambda_0) \subset
S$ containing $0 \in S$, i.e. there exists a holomorphic section
$\Lambda \in  \Gamma(S(\lambda_0),\cF^2)$ such that $\lambda_t=\Lambda|{X_t}
\in F^2 \cap H^2(U_t,\Q)$ for all $t \in S(\lambda_0)$. Let $T_0 \subset
T_{S,0}$ be the holomorphic tangent space to $S(\lambda_0)$ at $0 \in S$.
Since $\lambda_0$ is an integral class,
$$\nabla^2(\lambda_0) \in {\rm Hom}(T_{S(\lambda_0),0},(\cF^1/\cF^2)_0) $$
satisfies $\nabla^2(\lambda_0)(T_0)=0$. Therefore -
if $W_0 \subset W(log)$ denotes the subspace corresponding to $T_0$
under the logarithmic Kodaira-Spencer map - we  have that
$$\C \cdot \lambda_0 \otimes W_0 \mapsto 0 \in H^1(X,\Omega_X^1(logZ)) $$
i.e. $\C \cdot \lambda_0$ is contained in the left kernel of $\nabla^2$
with respect to $W_0$. By the assumption, however, this implies that $W_0$
is a proper subspace of $W(log)$ and hence also $T_0$ is a proper subspace
of $T_{S,0}$. It follows that $S(\lambda_0)$ is a proper subvariety of
$S$. The countability follows from the countability of the
groups $H^2(U_t,\Q)$. \qed

\subsection{Proof of Proposition 3.4.}

We use the notations and assumptions of the main theorem and let
$k$ be the number of intersection points of $Z_1$ and $Z_2$. \\

\prf (of Prop. 3.4.)\\
Let $X_0:=X \setminus Z_{sing}$ and $U=X \setminus Z$,
where $Z=Z_1 \cup  Z_2$.
By weak purity $H^i_{\cD,|Z_{sing}|}(X,\Z(2))=0$
for $i \le 3$. Therefore  $H^i_{\cD}(X,\Z(2)) =H^i_{\cD}(X_0,\Z(2))$ for
$i \le 2$ und $H^3_{\cD}(X,\Z(2)) \subset H^3_{\cD}(X_0,\Z(2))$. Note that
$H^4_{\cD,|Z_{sing}|}(X,\Z(2)) \cong \Z^k$. Also let
$K:={\rm Ker}(H^4_{\cD,|Z_{sing}|}(X,\Z(2)) \to H^4_{\cD}(X,\Z(2))) $ \\
There is a diagram
$$\matrix{
&& 0 && 0 &&\cr
&& \downarrow && \downarrow \cr
H^2_{\cD}(U,\Z(2)) & \to &H^3_{\cD,|Z|}(X,\Z(2)) &
\to & H^3_{\cD}(X,\Z(2)) &\to & H^3_{\cD}(U,\Z(2)) \cr
|| & & \downarrow & & \downarrow & & || \cr
H^2_{\cD}(U,\Z(2)) & \to &H^3_{\cD,|Z|}(X_0,\Z(2)) &
\to & H^3_{\cD}(X_0,\Z(2)) &\to & H^3_{\cD}(U,\Z(2)) \cr
&& \downarrow && \downarrow && \cr
&& K & = & K && }
$$
Let
$D_i=Z_i \setminus Z_{sing}$. The $D_i$ are smooth disjoint divisors in $X_0$
and one has an isomorphism $H^3_{\cD,|Z|}(X_0,\Z(2)) \cong
\oplus H^1_{\cD}(D_i,\Z(1))$. These groups can be computed via the
exact sequences
$$ 0 \to \C^* \to H^1_{\cD}(D_i,\Z(1)) \to F^1 \cap H^1(D_i,\Z) \to 0
$$
for $i=1,2$. The localization sequence
$$ 0 \to \C^*=H^1_\cD(Z_i,\Z(1)) \to H^1_{\cD}(D_i,\Z(1)) \to
H^0_\cD(Z_{\rm sing},\Z(0)) \cong \Z^k \to H^2_\cD(Z_i,\Z(1)) $$
identifies $F^1 \cap H^1(D_i,\Z)$ with ${\rm Ker}(\Z^k \to
{\rm Pic}(Z_i))$. Using the fact that \\
${\rm Ker}(\oplus_i {\rm Ker}(\Z^k \to {\rm Pic}(Z_i)) \to K) =
{\rm Ker}(\Z^k \to \oplus_i {\rm Pic}(Z_i) )$
we have proved the exactness of
$$0 \to (\C^*)^{\oplus 2} \to  H^3_{\cD,|Z|}(X,\Z(2)) \to {\rm Ker}(\Z^{k}
\to \oplus {\rm Pic}(Z_i))  \to 0 $$
and it remains to show that $CH^1(|Z|,1) \cong H^3_{\cD,|Z|}(X,\Z(2))$.
But the $D_i$ are smooth and hence $CH^1(D_i,1) \cong H^1_{\cD}(D_i,\Z(1))$
and the localization sequence for
$D_i \subset Z_i$ with complement  $Z_{sing}$ gives an exact sequence
$$0 \to (\C^*) \cong CH^1(Z_i,1) \to  CH^1(D_i,1) \to CH^0(Z_{\rm sing},0)
\cong \Z^{k} \to CH^1(Z_i) \to \ldots $$
and the claim follows in the same way as for $H^3_{\cD,|Z|}(X,\Z(2))$.  \qed

\subsection{Proof of the Main Theorem}

We use the notations as in the theorem.\\
{\bf Proof of the main theorem:} By the proposition above we have a sequence
$$0 \to (\C^*)^{\oplus 2} \to CH^1(|Z|,1) {\buildrel \tau \over \to}
{\rm Ker}(\Z^{k} \to \oplus {\rm Pic}(Z_i)) \to 0    $$
and assumption (1) is equivalent to $\tau(Z) \neq 0$. The complex
$$ H^2_{\cD}(X,\Z(2)) \to H^2_{\cD}(U,\Z(2)) \to H^3_{\cD,|Z|}(X,\Z(2))
\to H^3_{\cD}(X,\Z(2)) \to H^3_{\cD}(U,\Z(2)) $$
has the subcomplex
$$ H^1(X,\C^*) \to H^1(U,\C^*) \to (\C^*)^2 \to
{\rm NS(X)} \otimes \C^*  \to 0 $$
by assumption (2d). Here we have assumed that $Z_1$ and $Z_2$ generate
${\rm NS}(X) \otimes \Q$, because by assumption (2d) we can always choose
a general deformation of $X$ that satisfies this property and
also the other assumptions of
the theorem, since they are of generic nature. We get a diagram:
$$\matrix{
CH^2(U,2) & \to & CH^1(Z,1) & \to & CH^2(X,1) & \to & CH^2(U,1) \cr
\downarrow && \downarrow \tau && \downarrow && \downarrow \cr
{H^2_{\cD}(U,\Z(2))\over H^1(U,\C^*)} &
{\buildrel \alpha \over \to} & {\rm Ker}(\Z^{k} \to
\oplus {\rm Pic}(Z_i)) & \to &
{H^3_{\cD}(X,\Z(2)) \over {\rm NS}(X) \otimes \C^*}  & \to
& H^3_{\cD}(U,\Z(2))    }
$$
Since $H^2_{\cD}(U,\Z(2))/ H^1(U,\C^*)=F^2 \cap H^2(U,\Z)$,
it is enough to show that the image of
$\alpha: F^2 \cap H^2(U,\Z) \to {\rm Ker}(\Z^{k} \to
\oplus {\rm Pic}(Z_i)) $ does
not contain $\tau(Z)$. Note that the kernel is a free group, so it
will be enough to show that $F^2 \cap H^2(U_t,\Q)=0 $
for general $t$ by the following argument.\\
{\it Claim:} The cup-product map
$H^0(X,\Omega^2_X(logZ)) \otimes W(log) \to H^1(X,\Omega_X^1(logZ)) $
has no left kernel.\\
Taking this for granted we are finished with the proof
by Lemma (5.2).\\
Hence it remains to prove the claim:\\
Look at the commutative diagram that exists by assumption (2c):
$$\matrix{
H^0(X,\Omega^2_X(logZ)) \otimes H^1(X,T_X(-Z)) & \to &
H^1(X,\Omega^1_X)/\oplus H^0(Z_i,\cO_{Z_i}) \cr
\downarrow & & \bigcap \cr
H^0(X,\Omega^2_X(logZ)) \otimes W(log) & \to & H^1(X,\Omega^1_X(logZ)) }
$$
By assumption (2b) the upper map has no left kernel and therefore also the
lower map, since $H^1(X,\Omega^1_X)/\oplus H^0(Z_i,\cO_{Z_i}) \to
H^1(X,\Omega^1_X(logZ))$ is an injection.        \qed

\section{Examples}

\subsection{Sufficiently Ample General
Complete Intersections in Arbitrary Varieties}

The result that follows has several roots. First of all there is the
Noether-Lefschetz theorem stated in the introduction, which says that a
general hypersurface of degree $d \ge 4$ in $\P^3_{\C}$ has ${\rm Pic}(X)
\cong \Z$. Then there is the theorem of Green-Voisin (\cite{Gre}), saying
that a general hypersurface of degree $d \ge 6$ in $\P^4_{\C}$ has torsion
Abel-Jacobi invariants. In 1989 we started to generalize this
together with M. Green by using
a sophisticated version of residue representations of differential forms
via global sections of adjoint linear systems. But at the same time M. Nori
came up with his fantastic connectivity theorem from \cite{Nori}, which gave
a much shorter proof of the following result:
\begin{theorem} \cite{GM} \\
Let $(Y,\cO(1))$ be a smooth and projective polarized variety
of dimension $n+h$, $X \subset Y$ a
general complete intersection of dimension $n$ and
multidegree $(d_1,..,d_h)$ with $min(d_i)$ sufficiently large. Furthermore
assume that $0 \le p \le n-1$. Then:
$${\rm Image}( CH^p(X) \otimes \Q \to H^{2p}_\cD(X,\Q(p))) $$
$$={\rm Image}(CH^p(Y) \otimes \Q \to H^{2p}_\cD(X,\Q(p))) $$
modulo possibly the image of a certain subtorus of $J^p(Y)_\Q$, which
vanishes if the generalized Hodge conjecture holds.
\end{theorem}
As always {\sl general} means for all points in the moduli space outside
a countable set of proper analytic subvarieties, which is sometimes
also called {\sl very general}.
The generalization to higher Chow groups can be proved with the same method
(this was also observed by S.Bloch, M.Nori and C. Voisin and
is mentioned in \cite{Voi1}):
\begin{theorem}
Let $(Y,\cO(1))$ be a smooth and projective polarized variety
of dimension $n+h$, $X \subset Y$ a
general complete intersection of dimension $n$ and
multidegree $(d_1,..,d_h)$ with $min(d_i)$ sufficiently large. Furthermore
assume that $1 \le p \le n$. Then:
$${\rm Image}( CH^p(X,1) \otimes \Q \to H^{2p-1}_\cD(X,\Q(p))) $$
$$ \subset {\rm Image}(H^{2p-1}_\cD(Y,\Q(p))  \to H^{2p-1}_\cD(X,\Q(p))) $$
\end{theorem}
\prf
Let $S:=\prod_i^h \P(H^0(Y,\cO_Y(d_i)))$ and $f:B \to S$ the universal
complete intersection. If $g: T \to S$ is any smooth morphism, we write
also $g: B_T=B \times_S T \to B$ for the base change and $A_T=Y \times T$.
If $Z_s$ is a cycle on $X=X_s$
for $s \in S$, we can find a cycle $\cZ \in CH^p(B_T,1)$
for some smooth morphism
$g: T \to S$ with the property that for some $t \in g^{-1}(s)$ we have
$\cZ \cap X_t=Z_s$ and such that $g$ is smooth and finite.
By Nori's theorem \cite{Nori} we have for $1 \le p \le n$
that the restriction homomorphism
$$  i^*: H^{2p-1}_\cD(A_T,\Q(p)) \longrightarrow H^{2p-1}_\cD(B_T,\Q(p))$$
is an isomorphism. Therefore there is a cohomology class
$\alpha \in H^{2p-1}_\cD(A_T,\Q(p))$ such that $i^* \alpha =c_{p,1}(\cZ)$.
The commutative diagram of restriction maps
$$\matrix{ H^{2p-1}_\cD(A_T,\Q(p)) & \to & H^{2p-1}_\cD(B_T,\Q(p)) \cr
\downarrow &  & \downarrow \cr
H^{2p-1}_\cD(Y,\Q(p)) & \to & H^{2p-1}_\cD(X_s,\Q(p)) }
$$
shows that $\alpha|_{X_t}$ lies in the image of
$H^{2p-1}_\cD(Y,\Q(p))$ for every $t$ with $g(t)=s$
and therefore the theorem is proved.  \qed

\begin{corollary}
Let $X \subset \P^3$ be a general hypersurface of degree
$d \ge 5$. Then: The Chern class $CH^2(X,1) \otimes \Q \to
H^3_\cD(X,\Q(2)) $ has image isomorphic to  $H^3_\cD(\P^3,\Q(2)) \cong
\C/\Q(1)$.
\end{corollary}
\prf Let $B_T$ be any smooth base change of the universal
hyperplane section $f: B \to S$ with $S=\P(H^0(\P^3,\cO_{\P^3}(d)))$.
By  \cite{Par} we have
$H^3_\cD(\P^3 \times T,\Q(2)) \cong  H^{3}_\cD(B_T,\Q(2))$ for
$d \ge 5$. This is the required equality in \cite{GM}. \qed

\subsection{K3 Surfaces}

\subsubsection{A Deformation of the Fermat Quartic}

Here we give an explicit example satifying the assumptions of the
main theorem:\\
Consider the Fermat quartic surface
$$X=\{x_0^4+x_1^4-x_2^4-x_3^4=0\} \subset \P^3$$
and the linear forms
$$L_0=x_1-x_2, \quad L_1=x_0-x_3, \quad L_2=x_1+x_3, \quad L_3=x_0-x_2 $$
in $\P^3$. We look at the family of quartics
$$F_t(x_0:x_1:x_2:x_3)=x_0^4+x_1^4-x_2^4-x_3^4 + 2t L_0L_1L_2L_3 $$
and define $X_t:=\{ x \in \P^3 \mid F_t(x)=0 \}$.\\
The following line on $X$ is important:
$$G=\{x_1-x_2=x_0-x_3=0\}
$$
Planes containing it define a pencil of planes in $\P^3$:
$$H_\lambda=\{x_1-x_2-\lambda(x_0-x_3)=0\}
$$
The residual curves $E_{\lambda,t}$ to the intersection of $H_\lambda$ with
$X_t$
are elliptic curves and their equations can be computed as follows: The
coordinate transformation $(x_0:x_1:x_2:x_3) \mapsto ((x:y:z),\lambda)$,
where
$$x_0=x+z,x_1=y+\lambda z,x_2=y-\lambda z ,x_3=x-z$$
defines a rational map (blowing up along $G$) and we obtain the new
equations
$$F_t(x,y,z,\lambda)=8z(x^3+\lambda y^3+z^2(x+\lambda^3 y) +\lambda t z
L_2 L_3) =:8z F'(x,y,z,\lambda)$$
where $L_2(x,y,z,\lambda)=x+y+(\lambda-1)z,
L_3(x,y,z,\lambda)=x-y+(\lambda+1)z$
and hence for the strict transforms we get the following equations in
$\P^2 \times \P^1$:
$$X_t=\{(x,y,z,\lambda) \mid x^3+\lambda y^3+z^2(x+\lambda^3 y)
+\lambda t z L_2(x,y,z,\lambda) L_3(x,y,z,\lambda) =0 \}$$
defining an elliptic pencil with fibers $E_{\lambda,t}$ for every $t$.\\
The strict transform of $G$ is $G=\{z=0\}$ and $E_{\lambda,t} \cap G =
\{P_1=(a:1:0),P_2=(a\zeta:1:0),P_3=(a\zeta^2:1:0)\}$ independent of $t$ with
$a^3:=-\lambda$ and $\zeta=e^{2\pi \sqrt{-1}/3} $. To compute the tangent
lines $T_i$ of $E_{\lambda,t}$ at $P_i$ we compute the partial derivatives
$$ {\partial F'_t \over \partial x} = 3x^2 +z^2 + \lambda t z (L_2+L_3)$$
$$ {\partial F'_t \over \partial y} = 3\lambda y^2 +\lambda^3 z^2
+ \lambda t z (L_3-L_2)$$
$$ {\partial F'_t \over \partial z} = 2z(x+\lambda^3 y) + \lambda t
L_2(x,y,z,\lambda) L_3(x,y,z,\lambda) + \lambda t z
((\lambda-1)L_3+(\lambda+1)L_2) $$
Therefore we get the equations of the tangent lines
$$T_1=\{ 3a^2 x + 3 \lambda y +\lambda t (a^2-1) z=0 \} $$
$$T_2=\{ 3a^2 \zeta^2 x + 3 \lambda y +\lambda t (a^2 \zeta^2 -1) z =0 \} $$
$$T_3=\{ 3a^2 \zeta^4 x + 3 \lambda y +\lambda t (a^2 \zeta^4- 1) z =0 \}$$
with intersection $T_1 \cap T_2 \cap T_3 =\{(-\lambda t:t:3) \} =:
\{P_{\lambda,t} \} $. The projection from $P_{\lambda,t}$ onto $G=\{z=0\}$
defines a hyperelliptic morphism onto $G$ if and only if
$P_{\lambda,t} \in E_{\lambda,t}$. This is one algebraic condition
we will compute:
$$0= F'(P_{\lambda,t})=t \lambda (\lambda-1) F''_t(\lambda) $$
with $F''_t=\lambda(t^2-18t+9)+4t^2-6t-18$. Solving $F''_t(\lambda)=0$
we get
$$\lambda(t) = -{4t^2-6t-18 \over t^2-18t+9}  $$

\subsubsection{Non-Trivial Cycles}

If $\lambda=\lambda(t)$, the hyperelliptic map of the preceeding section
$E_{\lambda(t),t} \to G$ defines
a rational function $f_1$ with zero divisor $2(P_1)-2(P_2)$ as a cycle on
$E_{\lambda(t),t}$. There is always a rational function $f_2$ on $G$
with zero divisor $2(P_2)-2(P_1)$ as a quotient of 2 squares. Let
$Z:=E_{\lambda(t),t} \otimes f_{1,t} + G \otimes f_{2,t} \in CH^2(X_t,1)$.
Note that the cycle $Z_t$ is only defined over a covering of $\P^1$
since the choice of two of the 3 intersection points with an order is only
determined up to permutations. \\
However in order to verify the assumptions in the main theorem,
in particular (2c), we have to extend the base space of our deformation. Let
$N \subset \P H^0(\P^3,\cO_{\P^3}(4))$ be an irreducible
component of the Noether-Lefschetz locus of quartics containing a line and
containing the deformation of the Fermat quartic in section (6.2.1).
After some smooth base change $g: S \to N$ we get a family of lines
$G_t$, elliptic curves $E_t$ in $X_t$, rational functions $f_{1,t}$ on $E_t$,
$f_{2,t}$on $G_t$ and cycles $Z_t=E_t \otimes f_{1,t} + G_t \otimes f_{2,t}
\in CH^2(X_t,1)$ well-defined on $S$ that extend the cycles
$E_{\lambda(t),t}$ and $G=G_0$ above. $S$ can be chosen in such a way that
assumption (2a) of the criterion is satisfied.
The cycles even define classes also denoted
by $Z_t$ in $CH^1(|Z_t|,1)$ that map to $CH^2(X_t,1)$ under the natural map
$CH^1(|Z_t|,1) \to CH^2(X_t,1)$. \\
$W(log) \subset H^1(X,T_X(logZ))$ is a codimension one subspace,
since the deformations of the pair $(X_t,Z_t)$
that preserve the cycle in $CH^2(X_t,1)$ map generically
finite onto
the deformations of $X$ itself: on every generic elliptic fibration
only a finite number of fibers will have rationally equivalent intersection
points with a line, as we checked in the example above. Note that
it is sufficient to check such an assertion at one point of the moduli space.
Thus $W(log) \to W$ is an isomorphism.
$W$ is isomorphic to the tangent space to $N$
at the point $X$ by Noether-Lefschetz theory. The following diagram is
therefore exact and commutative:
$$\matrix{ 0 \to \C & \to & H^1(X,T_X(log Z)) & \to &  H^1(X,T_X) &
\to & \C^2 & \to & 0 \cr
&& \bigcup && \bigcup &&&& \cr
&& W(log) & {\buildrel \cong \over \to} & W &&&& }
$$
\begin{lemma} The cycle
$Z_t=E_{t} \otimes f_{1,t} + G_t \otimes f_{2,t}$ in
$CH^2(X_t,1)/{\rm Pic}(X_t)\otimes \C^*$ is nontrivial modulo torsion
for $t$ outside a countable union of proper analytic subsets of $S$.
\end{lemma}
\prf  We have to verify the assumptions of the main theorem.
To simplify the notation, we assume for the moment that $t=0$ and let
$E=E_0,G=G_0$ and $Z=Z_0$. $Z$ will also denote the effective cycle $E+G$.
To check (1), it is sufficient to remark that $\tau(Z)$ corresponds to the
element $(2,-2,0)$ in ${\rm Ker}(\Z^3 \to {\rm Pic}(E) \oplus {\rm Pic}(G))$
and therefore does not decompose. Note that if we write
$\tau(Z)$ as $(2,-2,0)$, we have chosen an ordering of our curves $Z_1$ and
$Z_2$ (as we did in the proof of prop. 3.4.), and then the map
$\tau: CH^1(|Z|,1) \to \Z^3$ is given by $Z_1 \otimes h_1+
Z_2 \otimes h_2 \mapsto {\rm div}(h_1)=-{\rm div}(h_2)$.\\
Assumption  2(a) is clear by construction.\\
To prove (2b), we look at the following graded Artinian ring $R_*=
\oplus R_d$ with
$$R_d={{H^0(\P^3,\cO_{\P^3}(d))} \over J_{F,d}}$$
In the example the Jacobian ideal $J_{F,d}$ is the d-th graded piece of the
homogenous ideal generated by the monomials $x_i^3$.
One has the following isomorphisms first proved by Griffiths (see
\cite{1594}, page 44):
$$H^1(X,\Omega_X^1)_{\rm pr} \cong R_4,\quad H^{0,2}(X) \cong R_8 \cong \C$$
$$H^1(X,T_X(-Z)) \cong R_3, R_0 \cong \C, R_1 \cong H^0(\Omega^2_X(logZ))
$$
Further there is Macaulay duality: $R_k^* \cong R_{8-k}$, see op.cit.\\
In this language we have to prove that
$$R_1 \otimes R_3 \longrightarrow R_4/\C \cdot g  $$
has no left kernel, where $g$ is a quartic polynomial representing the
extra cohomology class of $G$ (note that the cohomology class of $Z=G+E$ is
the hyperplane class and therefore not primitive). In other words we have to
prove that
$$R_1 \longrightarrow R_3^* \otimes R_4/ \C \cdot g $$
is injective.  After dualizing this means that
$$m: V/J_{F,4} \otimes R_3 \longrightarrow R_7 $$
is surjective, where $V \subset H^0(\P^3,\cO_{\P^3}(4))$ is the codimension
one linear subspace dual to $\C \cdot g$. $V$ contains $J_{F,4}$ and is
therefore basepointfree. The surjectivity of $m$ then follows from
the theorem of Green-Gotzmann (\cite{1594}, pg.74),
which has the following corollary:\\
Let $V \subset H^0(\P^3,\cO_{\P^3}(d))$ be a linear basepointfree subspace
of codimension $c$, then the map
$$m : V \otimes H^0(\P^3,\cO_{\P^3}(k-d)) \to H^0(\P^3,\cO_{\P^3}(k))$$
is surjective if $k \ge d+c$. \\
We apply this result for $k=7,d=4,c=1$ thus proving (2b).\\
To prove (2c), consider the sheaf $\cG:=\Omega^2_X(logZ)/\Omega^2_X$, which
we can also denote by $\omega_Z$. One has an exact sequence
$$0 \to \Omega^1_G \oplus \Omega^1_E \to \cG \to \C^{\oplus 3} \to 0$$
which, since $G$ is rational, has the same $H^0$ sequence of vector spaces as
the one coming from
$$0 \to \Omega^1_E \to \Omega^1_E(log G) \to \C^{\oplus 3}\to 0 $$
and the coboundary maps are equal.
Hence here  $H^0(\cG) \cong H^0(E,\Omega^1_E(log G))$.
The cycle
$Z \in CH^1(|Z|,1)$ defines a non-zero element $\alpha \in H^0(\cG)$ with
$\alpha={1 \over {2\pi i}} \cdot {df_1 \over f_1}$. There is an induced map
$$H^0(\cG) \otimes H^1(X,T_X(logZ)) \to H^1(E,\cO_E)  $$
(since $G$ is rational) and $W(log)$ equals the right annihilator
of $\alpha$ under that map (this follows, because $W(log) $ is a codimension
1 subspace of $H^1(X,T_X(logZ))$ that annihilates $\alpha$ and therefore is
equal to the annihilator, since the latter also has codimension 1). \\
To show (2c) it is sufficient to look at the
following diagram (note that the right column is exact but not the left one):
$$\matrix{
H^0(X,\Omega^2_X(logZ)) \otimes H^1(X,T_X(-Z)) & \to &
H^1(X,\Omega^1_X)/\C^2 \cr
\downarrow & & \downarrow \cr
H^0(X,\Omega^2_X(logZ)) \otimes H^1(X,T_X(logZ)) & \to &
H^1(X,\Omega^1_X(logZ)) \cr
\downarrow & & \downarrow \cr
H^0(\cG) \otimes H^1(X,T_X(logZ)) &\to &H^1(E,\cO_E)
}
$$
It shows that the image of
$H^0(X,\Omega^2_X(logZ)) \otimes H^1(X,T_X(-Z)) \to H^1(E,\cO_E)$ is zero,
hence the image of $H^1(X,T_X(-Z))$ in $H^1(X,T_X(logZ))$
annihilates $\alpha$ and therefore is a subspace of $W(log)$.
This finishes the proof of (2c). \\
To prove (2d), look at the multiplication map
$$W \otimes H^1(X,\Omega^1_X) \to H^2(X,\cO_X) \cong \C $$
$W$ annihilates the cycle classes of the curves $G$ and $E$ and the
restriction to the quotient
$$W \otimes {H^1(X,\Omega^1_X) \over {\C [G] \oplus \C [E]}}
\to \C $$
is a perfect pairing by the Hodge Riemann bilinear relations. It follows
that no class in $ {H^1(X,\Omega^1_X) \over {\C [G] \oplus \C [E]}} $
survives in a general deformation direction of $W$ by the same argument
as in the proof of lemma (5.2.). Therefore for a general
surface $X_t$ one has that $NS(X_t) \otimes \Q$ is generated by the two
elements $[G]$ and $[E]$.  \qed

\rmk With some more work, using a monodromy argument of H. Clemens which
was also used by A. Collino in \cite{Col}, one can probably prove infinite
generation for the image of $CH^2(X_t,1)$ for a general quartic hypersurface
containing a line. Since for oother examples
this was proved in \cite{Col} and \cite{Voi2}
we refrain from presenting it here. Obviously the idea would be to study the
monodromy around a countable set of loci on the parameter space
of the surfaces $X_t$.

\subsection{Examples of General Type}

In the previous section we have studied some special quartic surfaces and
have seen that on can apply the main theorem. Now we would like to show
that one can do a similar construction on some quintic hypersurfaces
$X$ in $\P^3$ and again apply the main theorem to get cycles that are
indecomposable in $CH^2(X,1)$. Note that they have to be very special
surfaces in view of the result (6.2).\\
To construct examples of general type, we
look at the Shioda hypersurface of degree 5 from \cite{Sh}:
$ X=\{x \in \P^3 \mid x_0x_1^4+x_1x_2^4+x_2x_0^4+x_3^5=0 \} $.
It has an automorphism $\sigma$ of order 65, given by
$(x_0:x_1:x_2:x_3) \mapsto (\zeta^{16}x_0:\zeta^{-4}x_1:\zeta x_2:x_3)$
where $\zeta$ is a 65-th root of unity.
Shioda proves that the Picard group of $X$ is of rank one, by showing
that $\sigma$ acts irreducibly on $H^2_{\rm pr}(X,\Q)$. \\
Let us look at the 1-parameter family
$$X_t=\{ x \in \P^3 \mid
F_t(x)=x_0 x_1^4+x_1 x_2^4 + x_2 x_0^4 + x_3^5 + tx_3 x_1^4 =0 \}
$$
for $t \in \P^1_\C$. Furthermore let $H_i=\{x_i=0\}$ be the coordinate
hyperplanes and define the curves $Z_{1,t}:=X_t \cap H_3$ and
$Z_{2,t}:=X_t \cap H_0$. Finally we set $P:=(0:1:0:0)$, $Q:=(0:0:1:0)$ and
$R:=(1:0:0:0)$.
\begin{proposition}

(a) $Z_{1,t}$ is smooth for all $t$, $Z_{2,t}$ is smooth for $t \neq 0,
\infty$ and $Z_{2,0}$ is a rational curve with a point of multiplicity 4 at
$P$.\\
(b) $P,Q \in Z_{1,t} \cap Z_{2,t}$.\\
(c) for all $t$: $4P=4Q$ in $CH_0(Z_{2,t})$.\\
(d) for all $t$: $13P=13Q$ in $CH_0(Z_{1,t})$.
\end{proposition}
\prf (a) We omit the subscripts t, if there is no ambiguity.
$Z_1 =\{ (x_0:x_1:x_2) \in \P^2 \mid
x_0 x_1^4+x_1 x_2^4 + x_2 x_0^4 =0 \}$ independent of $t$,
which is smooth by the Jacobian criterion. $Z_{2,t}$ depends on $t$:
$ Z_{2,t} =\{ (x_1:x_2:x_3) \in \P^2 \mid
x_1 x_2^4 + x_3^5 + t x_3 x_1^4 =0 \}$. The gradient is therefore
$ (x_2^4 + 4t x_3 x_1^3 : 4 x_2^3 x_1 : 5 x_3^4 + t x_1^4 )$, which is
nowhere zero for $t \neq 0,\infty$.
For $t=0$ one gets a 4-fold point at $P$.
$Z_{2,0}$ is rational by the parametrization $(t_0:t_1) \mapsto
(0: -t_0^5: t_1^5 : t_1^4 t_0) $. This proves (a).\\
(b) One computes
$Z_1 \cap Z_2 = X_t \cap H_0 \cap H_3 =
\{x \in \P^3 \mid x_0=x_3=0,\quad x_1 x_2^4=0 \}$. Hence
$Z_2 \cdot H_3 = Q + 4 P $ as cycles on $\P^3$ with multiplicity counted.
(b) follows.\\
(c) We obtain $Z_2 \cap H_1 = X_t \cap H_0 \cap H_1 =
 \{x \in \P^3 \mid x_0=x_1=0,\quad x_3^5=0 \}$. Hence
$Z_2 \cdot H_1 = 5Q $ as a cycle on $\P^3$. But $Z_2 \cdot H_3 $ is
rationally equivalent to $Z_2 \cdot H_1$ on $Z_2$. Together with (b), we
obtain that $Q+4P=5Q$ and hence $4P=4Q$ in  $CH_0(Z_{2,t})$ for all $t$.\\
(d) $Z_{1,t} \cap H_2 = X_t \cap H_3 \cap H_2 =
  \{x \in \P^3 \mid x_2=x_3=0,\quad x_0 x_1^4=0 \}$. It follows that
$Z_1 \cdot H_2 = P+4R$ as cycles. Finally we compute
$Z_1 \cdot H_1 = R+4 Q$ by symmetry. Together we get that
$Q+4P=P+4R=R+4Q$ in $CH_0(Z_1)$. Solving these equations gives that
$13P=13Q=13R$ in $CH_0(Z_1)$.\qed

\begin{corollary}
For all $t$ , $52P$ and
$52Q$ are rationally equivalent on both curves. This
defines rational functions $f_i$ on $Z_i$ with
${\rm div}(f_1)=52P-52Q=-{\rm div}(f_2)$ as a zero cycle on
$X_t$ for all $t$ and hence a cycle $Z_t \in
CH^2(X_t,1)$.
\end{corollary}

As in the previous example we now choose a maximal irreducible
component $N \subset
\P(H^0(\P^3,\cO(5)))$, containing the 1-parameter family above and
such that the properties of the cycles $Z_i$
extend along a suitable \'etale covering $S$ of $N$ and therefore
cycles $Z_t \in CH^2(X_t,1)$ are well defined for all $t \in S$.
The assumptions of the main theorem can be checked in the
same way as in the previous example:
\begin{theorem}
For $t \in S$ general, $Z_t$ is indecomposable in $CH^2(X_t,1)$.
\end{theorem}
{\bf Proof:} Again the assumptions (1),(2a) are easy to verify. To prove (2b)
we use the same method as above. Assume $t=0$ and look at the Shioda
hypersurface. Here $K_X=\cO_X(1)$ and $K_X+Z=\cO_X(3)$. Again let $R_*$
denote the Griffiths' Jacobian ring. One has $H^0(X,K_X+Z)=R_3$ and
$H^{1,1}_{\rm pr}=R_6$. One also computes that $H^1(X,T_X(-Z)) = R_3$.
Note that $Z_1$ and $Z_2$ have the same cohomology class. Therefore we have
to show that $R_3 \otimes R_3 \to R_6$ has no left kernel, or dually that
$R_6 \otimes R_3 \to R_9 $ is surjective, which is obvious. This proves
(2b). The proof of (2c) proceeds in the same way as above with the necessary
modifications.   Finally (2d) holds since it holds on the Shioda surface $X_0$
and therefore the Picard number is also one on a general surface parametrized
by $S$. \qed

\section{Further Remarks}

First let us imitate a method of
Bloch (\cite{Bl3}) to show what structure the groups that are involved
should have. Note that $CH^2(X,1) \cong H^1(X,\cK_2)$ via the
Gersten-Quillen resolution.
Bloch's method is to study $H^1(X,\cK_2)$ infinitesimally
by giving an infinitesimal version of the sheaf $\cK_2$. Instead of $\C$ let
the ground field be $k$ and look at the local $k-$ algebra of dual numbers
$k[\epsilon]/\epsilon^2$. Imitating Bloch's construction in the case of
$CH^2(X)$ , it makes sense to speak of the formal tangent space
$TH^1(X,\cK_2)$
to $H^1(X,\cK_2)$ by defining it as the group $H^1(X,\cK_{2,\epsilon})$ where
$\cK_{2,\epsilon}$ is this case given by the absolute differentials
$\Omega^1_{X/\Q}$. This is also nicely explained in Murre's lecture in
\cite{1594}, chapter V. From the exact sequence
$$ 0 \to \cO_X \otimes_k \Omega^1_{k/\Q} \to
\Omega^1_{X/\Q} \to \Omega^1_{X/k} \to 0
$$
we get a long exact sequence
$$H^1(X,\cO_X) \otimes_k \Omega^1_{k/\Q} \to TH^1(X,\cK_2) \to
H^1(X,\Omega^1_{X/k}) \to \ldots $$
$$ \ldots \to H^2(X,\cO_X) \otimes_k \Omega^1_{k/\Q} \to TCH^2(X) \to
T{\rm Alb}(X)$$

We see two principles, however without giving a correct proof:
If $k=\C$ and $p_g(X)=0$ or $k=\Q$, then the Albanese map should
be injective (Bloch's conjecture on $CH^2(X)$) and if $k=\C$ and
$X$ has maximal Picard number, then the map ${\rm Pic}(X) \otimes \C^* \to
CH^2(X,1)$ should be surjective modulo torsion. But there is a caveat:
In \cite{Ram} there is an example of a submotive of
a Hilbert modular surface that has maximal Picard number but
indecomposable $CH^2(X,1)$. Over
number fields the situation seems to be different, we refer to
\cite{Ram}, \cite{Mil} and \cite{Fla}. I was told that also \cite{Mil}
is probably a counterexample. See also the discussion in \cite{Ras}.
\\
Another hope would be to give a direct computation
of the Chern class map $c_{2,1}$ in terms of integrals and to show that
for some special value of $t$, the Deligne classes of $Z_t$ are not torsion
modulo ${\rm NS}(X_t) \otimes \C^*$.\\
{\bf Question:} Does a direct computation of the integral imply that already
on the Shioda surface $X$ the cycle $Z$ is indecomposable? Note that
$X$ is defined over $\Q$ and a positive answer would contribute to the study
of the conjectures of Bloch and Beilinson. \\
It would also be nice to prove directly (without use of
Deligne cohomology) that $CH^2(X,1)/{\rm Pic}(X) \otimes \C^*$ is
only countable, for example by Hilbert scheme methods.
This was suggested to me by C. Voisin.\\
Finally one would like to study higher Chow groups of degenerations of
algebraic varieties, for example the degeneration of a quartic K3 surface
into a tetrahedron of 4 planes as we described in section (6.2.1). We would
like to know whether already a general member of that one-parameter
degeneration has an indecomposable $CH^2(X_t,1)$. This can probably be
done by using another slight variant of our
main theorem and checking the assumptions on the variety $X_{\infty}$, by
showing that there are no two-forms of a certain logarithmic type.
More general this should lead to a good understanding of the mixed motive of
a variety like the tetrahedron (consisting out of linear varieties with
certain combinatorial data) in the sense of \cite{Ja1} and relating it
to the smooth case via deformation theory as developed in chapter 3.\\

{\bf Acknowledgements}\\
It is a pleasure to thank H\'el\`ene Esnault for her encouragement and
support during this project and for the opportunity to
include theorem 4.2. into this paper.  \\
Special thanks go to Spencer Bloch, Mark Green, Marc Levine, Wayne Raskind,
Claire Voisin and James Lewis for some discussions and proof reading
and to Alberto Collino for
showing me his paper \cite{Col} at an early stage.\\
Financial support during this project was provided by the DFG
Forschergruppe in Essen and the DFG Schwerpunkt Komplexe Mannigfaltigkeiten.
The University of Chicago and the University of Leiden have offered support
for visits in October 1993 and October 1994 .

{\tt University Essen, e-mail: mueller-stach at uni-essen.de}


\begin{thebibliography}{}

\bibitem{BMS} {\sc F. Bardelli~} and {\sc ~S. M\"uller-Stach}: Infinitesimal
invariants of extensions of Hodge structures, unpublished (1991).

\bibitem{Be1} {\sc A.A.Beilinson}: Higher regulators and values of
L-functions, {\it Jour. Sov. Math.} {\bf 30} (1985),2036--2070.

\bibitem{Be2} {\sc A.A.Beilinson}: Notes on absolute Hodge cohomology,
 {\it Cont. Math.} {\bf 10} (1986),1--34.

\bibitem{Bl1} {\sc S.Bloch}: Algebraic cycles and higher K-theory,
 {\it Adv. in Math.}{\bf 61} (1986), 267-304; the moving lemma appears in
 {\it Jour. of Alg. Geometry} {\bf 3} (1994),537--568.

\bibitem{Bl2} {\sc S.Bloch}: Algebraic cycles and the Beilinson conjectures,
 {\it Cont. Math.} {\bf 58} (1) (1986), 65--79.

\bibitem{Bl3} {\sc S.Bloch}: Lectures on algebraic cycles,{\it Duke Univ.
Math. series} {\bf 4}(1980).

\bibitem{BKL} {\sc S.Bloch},{\sc ~A.Kas} and {\sc D.Liebermann}: Zero cycles
on surfaces with $p_g=0$, {\it Comp. Math.} {\bf 33} (1976), 135--145.

\bibitem{BO} {\sc S.Bloch~}and{\sc~A.Ogus}: Gersten's
conjecture and the homology of schemes, {\it Ann.Sci.Ecole
Norm.Sup.} {\bf 7} (1974) 181--202.

\bibitem{BS} {\sc S.Bloch}~and{\sc~V.Srinivas}: Remarks on
correspondences and algebraic cycles, {\it Amer. J.
Math.} {\bf 105} (1983), 1235--1253.

\bibitem{Col} {\sc A.Collino}: Griffiths' infinitesimal invariant and
higher K-theory on hyperelliptic Jacobians, preprint(1995).

\bibitem{Del} {\sc P.Deligne}: Th\'eorie de Hodge
II,III , {\it Publ. Math.}\, IHES {\bf 40} (1972), 5--57,
and {\bf 44} (1974), 5--78.

\bibitem{DS} {\sc C.Deninger}~and {\sc~A. Scholl}: The Beilinson conjectures,
{\it in Proceedings Cambridge Math. Soc. (eds. Coates and Taylor)} {\bf 153}
(1992),173--209.

\bibitem{EL} {\sc H.Esnault~}and {\sc~M.Levine}, Surjectivity of cycle maps,
{\it Proceedings Orsay 1992}, Ast\'erisque {\bf 218} (1993), 203--226.

\bibitem{E1}{\sc H.Esnault}: A note on the cycle map, {\it
Crelle J.} {\bf 411} (1990), 51--65.

\bibitem{E2}{\sc H.Esnault}: Une remarque sur la
cohomologie du faisceau de Zariski de la $K$-th\'eorie de
Milnor sur une vari\'et\'e lisse complexe, {\it Math.Z.}
{\bf 205} (1990), 373--378.

\bibitem{ES} {\sc H.Esnault} and {\sc V. Srinivas}: Chern classes of vector
bundles with holomorphic connections on a complete smooth complex variety,
{\it Jour. Diff. Geo.} {\bf 36} (1992),257--267.

\bibitem{EV}{\sc H.Esnault}~and{\sc~E.Viehweg}:
Deligne-Beilinson cohomology, in  Beilinson's
conjectures on special values of $L$-functions, {\it
Academic Press Perspectives in Math.} {\bf 4}
(1988) 43--92.

\bibitem{Fla}{\sc M.Flach}: A finiteness theorem for the symmetric square
of an elliptic curve,{\it Invent. math.} {\bf 109} (1992), 307--327.

\bibitem{Gil} {\sc H.Gillet}: Deligne homology and the
Abel--Jacobi maps, {\it Bull.AMS} {\bf 10}
(1984) 285--288.

\bibitem{Gre} {\sc M. Green:} Griffiths' infinitesimal invariant and
Abel-Jacobi maps, {\it Jour. Diff. Geom.}{\bf 29}(1989), 545-555.

\bibitem{1594} {\sc  M.Green et al.}: Algebraic cycles and Hodge theory,
{\it Springer LNM} {\bf 1594} (1994).

\bibitem{GM} {\sc M.Green~} and {\sc~S. M\"uller-Stach}:
Algebraic cycles on a general complete intersection of high multi-degree of
a smooth projective variety, to appear in {\it Comp. Math.} (1995).

\bibitem{Ja1} {\sc U.Jannsen}: Mixed motives and algebraic
$K$-theory, {\it Springer} LNM {\bf 1400} (1990).

\bibitem{Ja2} {\sc U.Jannsen}: Deligne homology, Hodge-D-conjecture and
motives, {\it Perspectives in Math.} {\bf 4} (1988), 305--372.

\bibitem{Ja3} {\sc U.Jannsen}: Motives and Filtrations on Chow groups,
{\it AMS proceedings motives conference at Seattle} {\bf 1} (1994), 245--303.

\bibitem{Lev} {\sc M.Levine}: Localization on singular varieties, {\it
Invent. math.} {\bf 91} (1988),423--464.

\bibitem{Lev2} {\sc M. Levine}: Bloch's higher Chow groups revisited, {\it
Proc. Conf. K-theory Strasbourg 1992}, Ast\'erisque 226 (1994).

\bibitem{Lew} {\sc J. Lewis:} A Survey of the Hodge conjecture,
University of Montreal Press (1990).

\bibitem{Mil} {\sc S.Mildenhall:} Cycles in a product of elliptic curves
and a group analogous to the cycle class group, {\it Duke J.} {\bf 67} (1992),
387--406.

\bibitem{Mum}{\sc D.Mumford}: Rational equivalence of
$0$-cycles on surfaces, {\it J.Math.Kyoto Univ.} {\bf 9}
(1968),195--204.

\bibitem{MS}{\sc S. M\"uller-Stach:} $\C^*-$extensions of tori, higher
Chow groups and
applications to equivalence relations for algebraic cycles, {\it Journal
of K-theory} {\bf 9} (1995), 395-406.

\bibitem{Nori} {\sc M.Nori}:Algebraic cycles and Hodge theoretic
connectivity, {\it Invent. math.} {\bf 111(1)} (1993),349--373.

\bibitem{Par} {\sc K.Paranjape}: Cohomological and cycle theoretic
connectivity, {\it  Ann. of math.} {\bf 140} (1994),641--660.

\bibitem{Ram} {\sc D. Ramakrishnan:} Letter to H. Esnault, (May 1994).

\bibitem{Ras} {\sc W.Raskind}: Algebraic K-theory, \'etale cohomology and
torsion algebraic cycles, {\it Contemp.Math.} {\bf 83} (1989),311--341.

\bibitem{Sch} {\sc P.Schneider}: Introduction to the Beilinson
conjectures, {\it Perspectives in Math.} {\bf 4} (1988), 1--36.

\bibitem{Sh} {\sc T. Shioda}: Algebraic cycles on Hypersurfaces in $\P^N$,
{\it Adv. Studies in Pure Math.{\bf 10}, Sendai Proc.} (1987), 717-732.
\bibitem{Sou} {\sc Ch.Soul\'e}: R\'egulateurs, {\it Bourbaki s\'eminaire
1984/85, exp. 644,ast\'erisque 133-134} (1986), 237--253.

\bibitem{Sus} {\sc A.Suslin} and {\sc A.Nesterenko}: Homology of the general
linear group over  a local ring and Milnor K-theory, {\it Izv. Akad. Nauk }
{\bf 34} (1990), 121--145

\bibitem{Sus} {\sc A.Suslin} and {\sc V.Voevodsky}: Singular homology of
abstract algebraic varieties, preprint (1994).

\bibitem{SZ} {\sc J.Steenbrink} and {\sc S. Zucker}: Variations of mixed
Hodge structure, I, {\it Inv. math.} {\bf 80} (1985), 489--542.

\bibitem{Voi1} {\sc C.Voisin}: Variations of Hodge structures
and algebraic cycles, lecture at ICM Z\"urich, preprint (1994).

\bibitem{Voi2} {\sc Ch.Oliva}: unpublished but described in \\
{\sc C.Voisin}:
Remarks on zero cycles of self-products of varieties, preprint (1995).

\bibitem{Voi3} {\sc C.Voisin}: Sur l'application d'Abel-Jacobi des
vari\'et\'es de Calabi-Yau de dimension trois, {\it Ann. scient. \'Ec.
normale superieure} {\bf 27}(1994), 209--226.
\end{thebibliography}
\end{document}